\documentclass[12pt,twoside]{amsart}
\usepackage{amssymb}
\usepackage{amscd}
\usepackage[abbrev,alphabetic]{amsrefs}
\usepackage{hyperref}
\usepackage{xypic}
\usepackage{xcolor}

\title[On del Pezzo fibrations in positive characteristic]
{On del Pezzo fibrations in positive characteristic} 
\author{Fabio Bernasconi and Hiromu Tanaka} 
\address{Department of Mathematics, Department of Mathematics, University of Utah, Salt Lake City, UT 84112, USA} 
\email{fabio@math.utah.edu}

\address{Graduate School of Mathematical Sciences, The University of Tokyo,
3-8-1 Komaba, Meguro-ku, Tokyo 153-8914, JAPAN}
\email{tanaka@ms.u-tokyo.ac.jp}

\newcommand{\NE}[0]{{\operatorname{NE}}}
\newcommand{\red}[0]{{\operatorname{red}}}

\newcommand{\Proj}[0]{{\operatorname{Proj}}}
\newcommand{\Spec}[0]{{\operatorname{Spec}}}
\newcommand{\Hom}[0]{{\operatorname{Hom}}}
\newcommand{\Bs}[0]{{\operatorname{Bs}}}

\newcommand{\pdeg}[0]{{\operatorname{p-deg}}}
\newcommand{\sep}[0]{{\operatorname{sep}}}

\newtheorem{thm}{Theorem}[section]
\newtheorem{lem}[thm]{Lemma}
\newtheorem{cor}[thm]{Corollary}
\newtheorem{prop}[thm]{Proposition}

\theoremstyle{definition}
\newtheorem{ex}[thm]{Example}

\newtheorem{dfn}[thm]{Definition}

\newtheorem{rem}[thm]{Remark}
      
\newtheorem{nota}[thm]{Notation}         
\newtheorem{step}{Step}

\makeatletter
  
  \@addtoreset{equation}{section}
  \makeatother

\newcommand{\MO}{\mathcal{O}}
\newcommand{\R}{\mathbb{R}}
\newcommand{\Q}{\mathbb{Q}}
\newcommand{\Z}{\mathbb{Z}}

\newcommand{\p}{\mathfrak{p}}
\newcommand{\m}{\mathfrak{m}}

\begin{document}

\begin{abstract}
We establish 
two results on three-dimensional del Pezzo fibrations in positive characteristic. 
First, we give an explicit bound for torsion index of relatively torsion line bundles. 
Second, we  show the existence of purely inseparable sections with explicit bounded degree. 
To prove these results, we study log del Pezzo surfaces defined over imperfect fields. 
\end{abstract}
\maketitle

\tableofcontents

\section{Introduction}

The minimal model conjecture predicts that 
an arbitrary algebraic variety is  
birational to either a minimal model or a Mori fibre space $\pi \colon V \rightarrow B$. 
A distinguished property of Mori fibre spaces in characteristic zero is that 
any relative numerically trivial line bundle is automatically trivial (cf. \cite[Lemma 3.2.5]{KMM87}). 
In \cite[Theorem 1.4]{Tan}, the second author constructs counterexamples 
to the same statement in positive characteristic. 
More specifically, if the characteristic is two or three, 
then there exists a Mori fibre space $\pi \colon V \rightarrow B$ and a line bundle $L$ on $V$ such that 
$\dim V=3, \dim B=1, L \equiv_{\pi} 0,$ and $L \not\sim_{\pi} 0$. 
Then it is tempting to ask how bad the torsion indices can be.  
One of the main results of this paper is to give such an explicit upper bound 
of torsion indices for three-dimensional del Pezzo fibrations.

\begin{thm}[Theorem \ref{t-triv-lb-mfs-curve}] \label{i-triv-lb-mfs-curve}
Let $k$ be an algebraically closed field of characteristic $p>0$.
Let $\pi \colon V \rightarrow B$ be a projective $k$-morphism such that $\pi_* \mathcal{O}_V = \mathcal{O}_B$, where $V$ is a three-dimensional $\mathbb{Q}$-factorial 
normal quasi-projective variety over $k$ and $B$ is a smooth curve over $k$. 
Assume there exists an effective $\mathbb{Q}$-divisor $\Delta$ such that $(V, \Delta)$ is klt and $\pi \colon V \rightarrow B$ is a $(K_V+\Delta)$-Mori fibre space.
Let $L$ be a $\pi$-numerically trivial Cartier divisor on $V$.
Then the following hold.
\begin{enumerate}
\item If $p \geq 7$, then $L \sim_{\pi} 0$.
\item If $p \in \left\{ 3, 5 \right\}$, then $p^2L \sim_{\pi} 0$.
\item If $p =2$, then $16 L \sim_{\pi} 0$.
\end{enumerate}
\end{thm}

We also prove a theorem of Graber--Harris--Starr type for del Pezzo fibrations in positive characteristic. 

\begin{thm}[Theorem \ref{t-rc-3fold}] \label{intro-rc-3fold}
Let $k$ be an algebraically closed field of characteristic $p>0$. 
Let $\pi:V \to B$ be a projective $k$-morphism such that $\pi_*\MO_V=\MO_B$, 
$V$ is a normal three-dimensional variety over $k$, and $B$ is a smooth curve over $k$. 
Assume that there exists an effective $\Q$-divisor $\Delta$ such that 
$(V, \Delta)$ is klt and $-(K_V+\Delta)$ is $\pi$-nef and $\pi$-big. 
Then the following hold. 
\begin{enumerate}
\item There exists a curve $C$ on $V$ such that $C \to B$ is surjective and 
the following properties hold. 
\begin{enumerate}
\item If $p\geq 7$, then $C \to B$ is an isomorphism. 
\item If $p \in \{3, 5\}$, then $K(C)/K(B)$ is a purely inseparable extension of degree 
$\leq p$. 
\item If $p=2$, then $K(C)/K(B)$ is a purely inseparable extension of degree 
$\leq 4$. 
\end{enumerate}
\item If $B$ is a rational curve, then $V$ is rationally chain connected. 
\end{enumerate}
\end{thm}

Theorem \ref{intro-rc-3fold} can be considered as a generalisation 
of classical Tsen's theorem, i.e. the existence of sections on ruled surfaces. 
Tsen's theorem was used to establish the log minimal model program in characteristic $p>5$ \cite[Section 3.4]{BW17}. 
Also, Tsen's theorem was used to show that 
$H^i(X, W\MO_{X, \Q})=0$ for 
threefolds $X$ of Fano type in characteristic $p>5$ when $i>0$ (cf. \cite[Theorem 1.3]{GNT}).

The proofs of Theorem \ref{i-triv-lb-mfs-curve} and Theorem \ref{intro-rc-3fold}
are carried out by studying the generic fibre $X:=V \times_B \Spec\,K(B)$ of $\pi$, which is a surface of del Pezzo type defined over an imperfect field.
Roughly speaking, Theorem \ref{i-triv-lb-mfs-curve} and Theorem \ref{intro-rc-3fold} 
hold by the following two theorems.

\begin{thm}[Theorem \ref{t-klt-bdd-torsion}]\label{i-klt-bdd-torsion}
Let $k$ be a field of characteristic $p>0$. 
Let $X$ be a $k$-surface of del Pezzo type. 
Let $L$ be a numerically trivial Cartier divisor on $X$. 
Then the following hold. 
\begin{enumerate}
\item If $p \geq 7$, then $L \sim 0$. 
\item If $p \in \{3, 5\}$, then $pL \sim 0$. 
\item If $p=2$, then $4L \sim 0$. 
\end{enumerate}
\end{thm}

\begin{thm}[Theorem \ref{t-ex-rat-points-dP}] \label{i-ex-rat-points-dP}
Let $k$ be a $C_1$-field of characteristic $p>0$. Let $X$ be a $k$-surface of del Pezzo type such that $k=H^0(X, \mathcal{O}_X)$.
Then
\begin{enumerate}
\item If $p \geq 7$, then $X(k) \neq \emptyset$;
\item If $p \in \left\{ 3,5 \right\}$, then $X(k^{1/p}) \neq \emptyset$;
\item If $p =2$ , then $X(k^{1/4}) \neq \emptyset$.
\end{enumerate}
\end{thm}

\subsection{Sketch of the proof of Theorem \ref{i-klt-bdd-torsion}}\label{ss-intro1}

Let us overview some of the ideas used in the proof of Theorem \ref{i-klt-bdd-torsion}. 
By considering the minimal resolution and running a minimal model program, 
the problem is reduced to the case when $X$ is a regular surface of del Pezzo type 
which has a $K_X$-Mori fibre space structure $X \to B$. 
In particular, it holds that $\dim B=0$ or $\dim B=1$.

\subsubsection{The case when $\dim B=0$} \label{sss2-intro0}

Assume that $\dim B=0$. 
In this case, $X$ is a regular del Pezzo surface. 
We first classify $Y:=(X \times_k \overline k)_{\red}^N$ 
(Theorem \ref{i-classify-bc}). 
We then compare $X \times_k \overline k$ with $Y=(X \times_k \overline k)_{\red}^N$ 
(Theorem \ref{i-p2-bound}).

\begin{thm}[Theorem \ref{t-classify-bc}]\label{i-classify-bc}
Let $k$ be a field of characteristic $p>0$. 
Let $X$ be a projective normal surface over $k$ with canonical singularities
such that $k=H^0(X, \MO_X)$ and $-K_X$ is ample.
Then the normalisation $Y$ of $(X \times_k \overline{k})_{\red}$ satisfies one of the following properties. 
\begin{enumerate}
\item $X \times_k \overline k$ is normal. 
Moreover, $X \times_k \overline k$ has at worst canonical singularities. 
In particular, $Y \simeq X \times_k \overline k$  and $-K_Y$ is ample. 
\item $Y$ is isomorphic to a Hirzebruch surface, i.e. a $\mathbb P^1$-bundle over $\mathbb P^1$. 
\item $Y$ is isomorphic to a weighted projective surface $\mathbb P(1, 1, m)$ 
for some positive integer $m$. 
\end{enumerate} 
\end{thm}

\begin{thm}[cf. Theorem \ref{t-p2-bound}]\label{i-p2-bound}
Let $k$ be a field of characteristic $p>0$. 
Let $X$ be a projective normal surface over $k$ with canonical singularities such that $k=H^0(X, \MO_X)$ and $-K_X$ is ample. 
Let $Y$ be the normalisation of $(X \times_k \overline k)_{\red}$ and let 
\[
\mu: Y \to X \times_k \overline k
\]
be the induced morphism. 
\begin{enumerate}
\item If $p \geq 5$, then $\mu$ is an isomorphism and $Y$ has at worst canonical singularities. 
\item If $p=3$, then the absolute Frobenius morphism $F_{X \times_k \overline k}$ 
of $X \times_k \overline k$ factors through $\mu$: 
\[
F_{X \times_k \overline k}:X \times_k \overline k\to Y \xrightarrow{\mu} X \times_k \overline k.
\]
\item 
If $p=2$, then the second iterated absolute Frobenius morphism $F^2_{X \times_k \overline k}$ 
of $X \times_k \overline k$ factors through $\mu$: 
\[
F^2_{X \times_k \overline k}:X \times_k \overline k\to Y \xrightarrow{\mu} X \times_k \overline k.
\]
\end{enumerate}
\end{thm}

Note that Theorem \ref{i-classify-bc} shows that $Y=(X \times_k \overline k)_{\red}^N$ 
is a rational surface. 
In particular, any numerically trivial line bundle on $Y$ is trivial. 
By Theorem \ref{i-p2-bound}, 
if $L'$ denotes the pullback of $L$ to $X \times_k \overline k$, 
then it holds that $L'^4 \simeq \MO_{X \times_k \overline k}$ in the case (3). 
Then the flat base change theorem implies that also $L^4$ is trivial. 

We now discuss the proofs of Theorem \ref{i-classify-bc}
and Theorem \ref{i-p2-bound}. 
Roughly speaking, we apply Reid's idea (\cite[cf. the proof of Theorem 1.1]{Rei94})
to prove Theorem \ref{i-classify-bc} 
by combining with a rationality criterion (Lemma \ref{l-rationality}). 
As for Theorem \ref{i-p2-bound}, 
we use the notion of Frobenius length of geometric normality $\ell_F(X/k)$ 
introduced in \cite{Tan19} (cf. Definition \ref{d-lF}, Remark \ref{r-lF}). 
Roughly speaking, if $p=2$, then we can prove that $\ell_F(X/k) \leq 2$ 
by computing certain intersection numbers (cf. the proof of Proposition \ref{p-p2-bound}). 
Then general result on $\ell_F(X/k)$ (Remark \ref{r-lF}) implies (3) of Theorem \ref{i-p2-bound}.

\subsubsection{The case when $\dim B=1$}\label{sss2-intro1}

Assume that $\dim B=1$, i.e. $\pi:X \to B$ is a $K_X$-Mori fibre space to a curve $B$. 
Since $X$ is of del Pezzo type, we have that the extremal ray $R$ of $\overline{\text{NE}}(X)$ 
that is not corresponding to $\pi:X \to B$ is spanned by an integral curve $\Gamma$, i.e. $R=\R_{\geq 0}[\Gamma]$. 
In particular, $\Gamma \to B$ is a finite surjective morphism of curves. 
If $K_X \cdot \Gamma<0$, then the problem is reduced 
to the above case (\ref{sss2-intro0}) by contracting $\Gamma$. 
Even if $K_X \cdot \Gamma=0$, then we may contract $\Gamma$ and 
apply the same strategy. 
Hence, it is enough to treat the case when $K_X \cdot \Gamma >0$. 
Note that the numerically trivial Cartier divisor $L$ on $X$ descends to $B$, 
i.e. we have $L \sim \pi^*L_B$ for some Cartier divisor $L_B$ on $B$. 
Then, a key observation is that 
the extension degree $[K(\Gamma):K(B)]$ is at most five (Proposition \ref{p-cov-deg-bound}). 
For example, if $p>5$, then $\Gamma \to B$ is separable. 
Then the Hurwitz formula implies that $-K_B$ is ample, hence $L_B \sim 0$. 
If $K(\Gamma)/K(B)$ is purely inseparable of degree $p^e$, 
then it hold that $L_B^{p^e} \sim 0$, since $-K_{\Gamma^N}$ is ample. 
For the remaining case, i.e. $p=2$, $[K(\Gamma):K(B)]=4$, and $K(\Gamma)/K(B)$ is inseparable but not purely inseparable, 
we prove that $H^0(B, L_B^4) \neq 0$ by applying Galois descent for the separable closure of $K(\Gamma)/K(B)$ (cf. the proof of Proposition \ref{p-ess-klt-bdd-torsion}).

\subsection{Sketch of the proof of Theorem \ref{i-ex-rat-points-dP}}\label{ss-intro2}

Let us overview some of the ideas used in the proof of Theorem \ref{i-ex-rat-points-dP}. 
The first step is the same as Subsection \ref{ss-intro1}, 
i.e. considering the minimal resolution and running a minimal model program, 
we reduce the problem to the case when $X$ is a regular surface of del Pezzo type 
which has a $K_X$-Mori fibre space structure $X \to B$. 

\subsubsection{The case when $\dim B=0$} 

Assume that $\dim B=0$. 
In this case, $X$ is a regular del Pezzo surface with $\rho(X)=1$. 
Since the $p$-degree of a $C_1$-field is at most one (Lemma \ref{l-Cr-pdeg}),
it follows from \cite[Theorem 14.1]{FS18} that $X$ is geometrically normal. 
Then Theorem \ref{i-classify-bc} implies that 
the base change $X \times_k \overline k$ is a canonical del Pezzo surface, 
i.e. $X \times_k \overline k$ has at worst canonical singularities 
and $-K_{X \times_k \overline k}$ is ample. 
In particular, we have that $1 \leq K_X^2 \leq 9$. 
Note that if $X$ is smooth, then it is known that $X$ has a $k$-rational point 
(cf. \cite[Theorem IV.6.8]{Kol96}). 
Following the same strategy as in \cite[Theorem IV.6.8]{Kol96}, 
we can show that $X(k) \neq \emptyset$ if $K_X^2 \leq 4$ 
(Lemma \ref{l-rat-pts-low-deg}).  
For the remaining cases $5 \leq K_X^2 \leq 9$, 
we use results established in \cite{Sch08}, 
which restrict the possibilities for the type of singularities on $X \times_k \overline k$. 
For instance, if $p \geq 11$, then \cite[Theorem 6.1]{Sch08} shows that 
the singularities on $X \times_k \overline k$ are of type $A_{p^e-1}$. 
However, such singularities cannot appear, because the minimal resolution 
$V$ of $X \times_k \overline k$ satisfies $\rho(V) \leq 9$. 
Hence, $X$ is actually smooth if $p \geq 11$ (Proposition \ref{p-dP-large-p1}). 
For the remaining cases $p \leq 7$, we study the possibilities one by one, 
so that we are able to deduce what we desire. 
For more details, see Subsection \ref{ss1-pi-pts}. 

\subsubsection{The case when $\dim B=1$} 

Assume that $\dim B=1$, i.e. $\pi:X \to B$ 
is a $K_X$-Mori fibre space to a curve $B$. 
Then the outline is similar to the one in (\ref{sss2-intro1}). 
Let us use the same notation as in (\ref{sss2-intro1}). 
The typical case is that $-K_B$ is ample. 
In this case, $B$ has a rational point. 
Then also the fibre of $\pi$ over a rational point, which is a conic curve, has a rational point. 
Although we need to overcome some technical difficulties, 
we may apply this strategy up to suitable purely inseparable covers 
for almost all the cases 
(cf. the proof of Proposition \ref{p-rat-point-mfs}). 
There is one case we can not apply this strategy: 
$p=2$, $K_X \cdot \Gamma>0$, and $K(\Gamma)/K(B)$ 
is inseparable and not purely inseparable. 
In this case, we can prove that $-K_B$ is actually ample (Proposition \ref{p-weird}).

\subsection{Large characteristic}

Using the techniques developed in this paper, we also prove the following theorem, which shows that some a priori possible pathologies of log del Pezzo surfaces over imperfect fields can appear exclusively in small characteristic.

\begin{thm}[cf. Corollary \ref{c-geom-red-7} and Theorem \ref{t-h1-vanish}]\label{intro-vanishing} 
Let $k$ be a field of characteristic $p \geq 7$. 
Let $X$ be a $k$-surface of del Pezzo type such that $k=H^0(X, \mathcal{O}_X)$. 
Then $X$ is geometrically integral over $k$ and $H^i(X, \MO_X)=0$ for any $i>0$. 
\end{thm}

As a consequence, we deduce the following result on del Pezzo fibrations in large characteristic:

\begin{cor} \label{c-genfb-largep}
Let $k$ be an algebraically closed field of characteristic $p \geq 7$. 
Let $\pi \colon V \to B$ be a projective $k$-morphism of normal $k$-varieties 
such that $\pi_* \MO_V= \MO_B$ and $\dim\,V-\dim\,B=2$.
Assume that there exists an effective $\Q$-divisor $\Delta$ on $V$ 
such that $(V,\Delta)$ is klt and $-(K_V+\Delta)$ is $\pi$-nef and $\pi$-big.
Then general fibres of $\pi$ are integral schemes and 
there is a non-empty open subset $B'$ of $B$ 
such that the equation 
$(R^i\pi_* \MO_V)|_{B'}=0$ holds for any $i>0$. 
\end{cor}

The authors do not know whether surfaces of del Pezzo type are geometrically normal if the characteristic is sufficiently large. 
On the other hand, even if $p$ is sufficiently large, 
regular surfaces of del Pezzo type can be non-smooth. 
More specifically, for an arbitrary imperfect field $k$ of characteristic $p>0$, 
we construct a regular surface of del Pezzo type which is not smooth 
(Proposition \ref{p-count}).

\subsection{Related results}

In this subsection, we summarise known results on log del Pezzo surfaces mainly over imperfect fields.

\subsubsection{Vanishing theorems} 

We first summarise results over algebraically closed fields of characteristic $p>0$. 
It is well known that smooth rational surfaces satisfy the Kodaira vanishing theorem (cf. \cite[Proposition 3.2]{Muk13}). 
However, the Kawamata--Viehweg vanishing theorem fails even 
for smooth rational surfaces (cf. \cite[Theorem 3.1]{CT18}). 
Moreover, the surface used in \cite[Theorem 3.1]{CT18} 
is a weak del Pezzo surface 
if the base field is of characteristic two (\cite[Lemma 2.4]{CT18}). 
Also in characteristic three, 
there exists a surface of del Pezzo type which violates the Kawamata--Viehweg 
vanishing (\cite[Theorem 1.1]{Ber}). 
On the other hand, if the characteristic is sufficiently large, 
it is known that surfaces of del Pezzo type satisfy the Kawamata--Viehweg vanishing by \cite[Theorem 1.2]{CTW17}. 

We now overview known results over imperfect fields. 
If the characteristic is two or three, 
there exists a surface $X$ of del Pezzo type such that $H^1(X, \MO_X) \neq 0$ 
(cf. Subsection \ref{ss1-patho}). 
On the other hand, regular del Pezzo surfaces of characteristic $p \geq 5$ 
satisfy the Kawamata--Viehweg vanishing theorem as shown in \cite[Theorem 1.1]{Das}.

\subsubsection{Geometric properties}

In characteristic two and three, 
there exist regular del Pezzo surfaces 
which are not geometrically reduced (cf. Subsection \ref{ss1-patho}). 
On the other hand, 
Patakfalvi and Waldron prove 
that regular del Pezzo surfaces are geometrically normal if the base field is of characteristic $p\geq 5$ (cf. \cite[Theorem 1.5]{PW}). 
Furthermore, 
Fanelli and Schr\"{o}er show that 
a regular del Pezzo surface $X$ is geometrically normal in every characteristic $p$
if $[k:k^p] \leq p$ and $\rho(X)=1$ (cf. \cite[Theorem 14.1]{FS18}).

\medskip

\indent \textbf{Acknowledgements:} 
We would like to thank P. Cascini, S. Ejiri, A. Fanelli, S. Schr\"{o}er, and J. Waldron for many useful discussions. 
We also would like to thank the referee for the constructive suggestions and reading the manuscript carefully. 
The first author was supported by the Engineering and Physical Sciences Research Council [EP/L015234/1].
The second author was funded by 
the Grant-in-Aid for Scientific Research (KAKENHI No. 18K13386). 

\section{Preliminaries}

\subsection{Notation}\label{ss-notation}

In this subsection, we summarise notation we will use in this paper. 

\begin{enumerate}
\item We will freely use the notation and terminology in \cite{Har77} 
and \cite{Kol13}. 
\item 
We say that a noetherian scheme $X$ is {\em excellent} (resp. {\em regular}) 
if 
the local ring $\MO_{X, x}$ at any point $x \in X$ is excellent (resp. regular). 
For the definition of excellent local rings, 
we refer to \cite[\S 32]{Mat89}. 
\item 
For a scheme $X$, its {\em reduced structure} $X_{\red}$ 
is the reduced closed subscheme of $X$ such that the induced morphism 
$X_{\red} \to X$ is surjective. 
\item For an integral scheme $X$, 
we define the {\em function field} $K(X)$ of $X$ 
to be $\MO_{X, \xi}$ for the generic point $\xi$ of $X$. 
\item 
For a field $k$, 
we say that $X$ is a {\em variety over} $k$ or a $k$-{\em variety} if 
$X$ is an integral scheme that is separated and of finite type over $k$. 
We say that $X$ is a {\em curve} over $k$ or a $k$-{\em curve} 
(resp. a {\em surface} over $k$ or a $k$-{\em surface}, 
resp. a {\em threefold} over $k$) 
if $X$ is a $k$-variety of dimension one (resp. two, resp. three). 
\item For a field $k$, we denote $\overline k$ (resp. $k^{\text{sep}}$) an algebraic closure (resp. a separable closure) of $k$. 
If $k$ is of characteristic $p>0$, 
then we set $k^{1/p^{\infty}}:=\bigcup_{e=0}^{\infty} k^{1/p^e}
=\bigcup_{e=0}^{\infty} \{x \in \overline k\,|\, x^{p^e} \in k\}$. 
\item For an $\mathbb{F}_p$-scheme $X$ we denote by $F_X \colon X \to X$ the {\em absolute Frobenius morphism}. For a positive integer $e$ we denote by $F^e_X \colon X \to X$ 
the $e$-th iterated absolute Frobenius morphism.
\item 
If $k$ is a field of characteristic $p>0$ such that $[k:k^p]<\infty$, we define its $p$-{\em degree} $\pdeg(k)$ as the non-negative integer $n$ such that $[k:k^p]=p^n$. 
The $p$-degree $\pdeg(k)$ is also called the degree of imperfection in some literature. 
\item If $k \subset k'$ is a field extension and $X$ is a $k$-scheme, we denote 
$X \times_{\Spec\,k} \Spec\,k'$ by $X \times_k k'$ or $X_{k'}$. 
\item 
Let $k$ be a field, let $X$ be a scheme over $k$ and 
let $k \subset k'$ be a field extension. 
We denote by $X(k')$ the set of the $k$-morphisms $\Hom_k(\Spec\,k', X)$. 
Note that if $X$ is a scheme of finite type over $k$ and $k \subset k'$ is 
a purely inseparable extension, then 
the induced map $\theta:X(k') \to X$ is injective and its image $\theta(X(k'))$ consists of closed points of $X$.  
\item Let $L$ be a Cartier divisor on a variety $X$ over $k$.
We define the {\em base locus} $\Bs(L)$ of $L$ 
by 
\[
\Bs(L):=\bigcap_{s \in H^0(X,L)} \left\{ x \in X \mid s(x)=0 \right\}.
\] 
In particular, $\Bs(L)$ is a closed subset of $X$. 
\item 
Let $k$ be an algebraically closed field. 
For a normal surface $X$ over $k$ and a canonical singularity $x \in X$ (i.e. a rational double point), we refer to the table at \cite[pages 15-17]{Art77} for the list of equations of type $A_n$, $D_n^m$ and $E_n^m$. 
For example, 
we say that $x$ is a canonical singularity of type $A_n$ 
if the henselisation of $\MO_{X, x}$ is isomorphic to $k\{x, y, z\}/(z^{n+1}+xy)$, 
where $k\{x, y, z\}$ denotes the henselisation of the local ring of $k[x, y, z]$ at 
the maximal ideal $(x, y, z)$. 
\end{enumerate}

\subsection{Geometrically klt singularities}
The purpose of this subsection is to introduce the notion of geometrically klt singularities and its variants. 

\begin{dfn}
Let $(X, \Delta)$ be a log pair over a field $k$ such that $k$ is algebraically closed in $K(X)$. 
We say that $(X, \Delta)$ is \emph{geometrically klt} (resp. terminal, canonical, lc) if 
$(X \times_k {\overline{k}}, \Delta\times_k{\overline{k}})$ is klt (resp. terminal, canonical, lc).
\end{dfn}

\begin{lem}\label{l-gred-open}
Let $k$ be a field. 
Let $X$ and $Y$ be varieties over $k$ which are birational to each other. 
Then $X$ is geometrically reduced over $k$ if and only if $Y$ is geometrically reduced over $k$. 
\end{lem}

\begin{proof}
Recall that for a $k$-scheme, 
being geometrically reduced is equivalent to being $S_1$ and geometrically $R_0$. 
Since both $X$ and $Y$ are $S_1$, 
the assertion follows from the fact that being geometrically $R_0$ is a condition on the generic point. 
\end{proof}

We prove a descent result for such singularities. 

\begin{prop}\label{p-klt-descent}
Let $(X, \Delta)$ be a geometrically klt (resp. terminal, canonical, lc) pair 
such that $k$ is algebraically closed in $K(X)$. 
Then $(X, \Delta)$ is klt (resp. terminal, canonical, lc).
\end{prop}

\begin{proof}
We only treat the klt case, as the others are analogous. 
Let $\pi \colon Y \rightarrow X$ be a birational $k$-morphism, where $Y$ is a normal variety and we write $K_Y+ \Delta_Y =\pi^*(K_X+\Delta)$. 
It suffices to prove that $\lfloor{\Delta_Y \rfloor} \leq 0$. 
Thanks to Lemma \ref{l-gred-open}, $Y$ is geometrically integral. 
Let $\nu \colon W \to Y \times_k \overline {k}$ be the normalisation morphism and let us consider the following commutative diagram:
\[
\begin{CD}
W  \\
@V \nu VV \\
Y \times_k \overline{k} @> g >> Y \\
@V \pi_{\overline{k}} VV @V \pi VV\\
X \times_k \overline{k} @> f >> X.
\end{CD}
\]
Denote by $\psi:= \pi_{\overline{k}} \circ \nu$ and $h:= g \circ \nu$ the composite morphisms.
We have
\[K_W+\Delta_W := \psi^*(K_{X_{\overline{k}}}+\Delta_{\overline{k}})=h^* \pi^*(K_X+\Delta)=h^* (K_Y + \Delta_Y). \]
By \cite[Theorem 4.2]{Tan18b}, there exists an effective $\mathbb{Z}$-divisor $D$ such that
\[h^*(K_Y+\Delta_Y)= K_Y+D+h^*\Delta_Y, \]
and thus $\Delta_W=D+h^*\Delta_Y \geq h^*\Delta_Y$. 
Since $(X_{\overline{k}}, \Delta_{\overline{k}})$ is klt, 
any coefficient of $\Delta_W$ is $<1$. 
Then any coefficient of $\Delta_Y$ is $<1$, thus $(X, \Delta)$ is klt.
\end{proof}

\begin{rem}
If $k$ is a perfect field, being klt is equivalent to being geometrically klt by
\cite[Proposition 2.15]{Kol13}. 
However, over imperfect fields, being geometrically klt is a strictly stronger condition. 
As an example, let $k$ be an imperfect field of characteristic $p>0$ and consider the log pair $(\mathbb{A}^1_k, \frac{2}{3}P)$, where $P$ is a closed point whose residue field $\kappa(P)$ is a purely inseparable extension of $k$ of degree $p$.
This pair is klt over $k$, but it is not geometrically lc.
\end{rem}

\subsection{Surfaces of del Pezzo type}

In this subsection, we summarise some basic properties of surfaces of del Pezzo type over arbitrary fields. 
For later use, we introduce some terminology. 
Note that del Pezzo surfaces in our notation allow singularities. 

\begin{dfn}\label{d-dP-wdP}
Let $k$ be a field. 
A $k$-surface $X$ is {\em del Pezzo} if 
$X$ is a projective normal surface such that $-K_X$ is an ample $\Q$-Cartier divisor. 
A $k$-surface $X$ is {\em weak del Pezzo} if 
$X$ is a projective normal surface such that $-K_X$ is a nef and big $\Q$-Cartier divisor. 
\end{dfn}

\begin{dfn}
Let $k$ be a field. 
A $k$-surface $X$ is \emph{of del Pezzo type} 
if $X$ is a projective normal surface over $k$ and 
there exists an effective $\mathbb{Q}$-divisor $\Delta \geq 0$ such that $(X, \Delta)$ is klt and $-(K_X+\Delta)$ is ample. 
In this case, we say that $(X, \Delta)$ is a log del Pezzo pair.
\end{dfn}

We study how the property of being of del Pezzo type behaves under birational transformations.

\begin{lem}\label{l-dP-min-res}
Let $k$ be a field. 
Let $X$ be a $k$-surface of del Pezzo type. 
Let $f : Y \to X$ be the minimal resolution of $X$.
Then $Y$ is a $k$-surface of del Pezzo type.
\end{lem}

\begin{proof}
Let $\Delta$ be an effective $\Q$-divisor such that $(X, \Delta)$ is a log del Pezzo pair.
We define a $\Q$-divisor $\Delta_Y$ by $K_Y+ \Delta_Y =f^*(K_X+\Delta)$.
Since $f:Y \to X$ is the minimal resolution of $X$, 
we have that $\Delta_Y$ is an effective $\Q$-divisor. 
The pair $(Y, \Delta_Y)$ is klt and $-(K_Y+ \Delta_Y)$ is nef and big. 
By perturbing the coefficients of $\Delta_Y$, 
we can find an effective $\Q$-divisor $\Gamma$ such that $(Y, \Gamma)$ is klt and $-(K_Y+\Gamma)$ is ample. 
\end{proof}

\begin{lem} \label{l-pert-ample}
Let $k$ be a field. 
Let $(X, \Delta)$ be a two-dimensional projective klt pair over $k$.
Let $H$ be a nef and big $\Q$-Cartier $\Q$-divisor.
Then there exists an effective $\Q$-Cartier $\Q$-divisor $A$ such that $A \sim_{\mathbb{Q}} H$ and 
$(X, \Delta + A)$ is klt.
\end{lem}

\begin{proof}
Thanks to the existence of log resolutions for excellent surfaces \cite{Lip78}, 
the same proof of \cite[Lemma 2.8]{GNT} works in our setting.
\end{proof}

\begin{lem}\label{l-dP-under-bir-mor}
Let $k$ be a field. 
Let $X$ be a $k$-surface of del Pezzo type. 
Let $f : X \rightarrow Y$ be a birational $k$-morphism to a projective normal $k$-surface $Y$. 
Then $Y$ is a $k$-surface of del Pezzo type.
\end{lem}
\begin{proof}
Let $\Delta$ be an effective $\Q$-divisor such that $(X, \Delta)$ is a log del Pezzo pair.
Set  $H:=-(K_X+\Delta)$, which is an ample $\Q$-Cartier $\mathbb{Q}$-divisor on $X$.  
By Lemma \ref{l-pert-ample}, there exists 
an effective $\Q$-Cartier $\Q$-divisor $A$ 
such that $A \sim_{\mathbb{Q}} H$ and  $(X, \Delta +A)$ is klt. 
Then the pair $(Y, f_*\Delta+f_*A)$ is klt and 
$K_X+\Delta+A \sim_{\mathbb{Q}} f^*(K_Y+f_*\Delta+f_*A) \sim_{\mathbb
Q} 0$. 
It follows from \cite[Corollary 4.11]{Tan18a} that $Y$ is $\Q$-factorial. 
By Nakai's criterion, the $\Q$-divisor $f_*A$ is ample. 
In particular $(Y, f_*\Delta)$ is a log del Pezzo pair.
\end{proof}

\subsection{Geometrically canonical del Pezzo surfaces}
In this subsection we collect results on the anti-canonical systems of geometrically canonical del Pezzo surfaces we will need later.

\subsubsection{Canonical del Pezzo surfaces over algebraically closed fields}

We verify that the results in \cite[Chapter III, Section 3]{Kol96} hold for del Pezzo surfaces with canonical singularities over algebraically closed fields. 
Recall that we say that $X$ is a canonical (weak) del Pezzo surface over a field $k$ if 
$X$ is a surface over $k$, $X$ is (weak) del Pezzo in the sense of Definition \ref{d-dP-wdP}, 
and $(X, 0)$ is canonical in the sense of \cite[Definition 2.8]{Kol13}.

\begin{prop} \label{p-cohomology-can-dP} 
Let $X$ be a canonical weak del Pezzo surface over an algebraically closed field $k$.
Then the following hold. 
\begin{enumerate}
\item $H^2(X, \MO_X(-mK_X))=0$ for any non-negative integer $m$. 
\item $H^i(X, \mathcal{O}_X) =0$ for any $i >0$. 
\item $H^0(X, \mathcal{O}_X(-K_X)) \neq 0$. 
\item $H^1(X ,\mathcal{O}_X(mK_X))=0$ for any integer $m$.  
\item $h^0(X, \mathcal{O}_X(-mK_X)) = 1+\frac{m(m+1)}{2} K_X^2$ 
for any non-negative integer $m$.  
\end{enumerate}
\end{prop} 

\begin{proof}
The assertion (1) follows from Serre duality. 
We now show (2). 
It follows from \cite[Theorem 5.4 and Remark 5.5]{Tan14} that 
$X$ has at worst rational singularities. 
Then the assertion (2) follows from the fact that $X$ is a rational surface 
\cite[Theorem 3.5]{Tan15}. 

We now show (3). 
By $H^2(X, \MO_X(-K_X))=0$ and the Riemann--Roch theorem, 
we have $h^0(X, \MO_X(-K_X)) \geq 1 + K_X^2 >0$.
Thus (3) holds. 

We now show (4).
By (3), there exists an effective Cartier divisor $D$ such that $D \sim -K_X$. 
In particular, $D$ is effective, nef, and big. 
It follows from \cite[Proposition 3.3]{CT} that 
\[
H^1(X, \MO_X(-nD))= H^1(X, \MO_X(K_X+nD))=0
\]
for any $n \in \Z_{>0}$. Replacing $D$ by $-K_X$, the assertion (4) holds. 
Thanks to (1) and (4), assertion (5) follows from the Riemann--Roch theorem. 
\end{proof}

\begin{lem}\label{l-comp-antican-sys}
Let $Y$ be a canonical weak del Pezzo surface over an algebraically closed field $k$. 
If a divisor $\sum_{i=1}^r a_i C_i \in |-K_Y|$ is not irreducible or not reduced, then every $C_i$ is a smooth rational curve.
\end{lem}

\begin{proof}
Taking the minimal resolution of $Y$, 
we may assume that $Y$ is smooth. 
Fix an index $1 \leq i_0 \leq r$. 
By adjunction, we have
\begin{equation}\label{e-comp-antican-sys}
2p_a(C_{i_0})-2 = -C_{i_0} \cdot \left(\sum_{i \neq i_0} \frac{a_i}{a_{i_0}} C_i \right) - \frac{a_{i_0}-1}{a_{i_0}} C_{i_0} \cdot (-K_Y).
\end{equation}
Note that both the terms on the right hand side are non-positive. 

Since $Y$ is smooth and $\sum_i a_i C_i$ is nef and big, 
it follows from \cite[Theorem 2.6]{Tan15} that $H^1(X, -n\sum_i a_i C_i)=0$ for $n \gg 0$. 
Hence, $\sum_i a_i C_i$ is connected. 
Therefore, if $\sum_i a_i C_i$ is reducible, 
the first term in the right hand side of (\ref{e-comp-antican-sys}) 
is strictly negative, hence $p_a(C_{i_0})<0$. 

If $a_{i_0} \geq 2$ and $C_{i_0} \cdot K_Y<0$, then 
the second term in the right hand side of (\ref{e-comp-antican-sys}) 
is strictly negative, hence $p_a(C_{i_0})<0$. 
If $C_{i_0} \cdot K_Y=0$, then 
$C_i$ is a smooth rational curve with $C_i^2=-2$. 
\end{proof}

\begin{prop} \label{p-gen-memb-antican}
Let $Y$ be a canonical weak del Pezzo surface over an algebraically closed field $k$. 
Let $\Bs(-K_Y)$ be the base locus of $-K_Y$, which is a closed subset of $Y$. 
Then the following hold. 
\begin{enumerate}
\item $\Bs(-K_Y)$ is empty or $\dim (\Bs(-K_Y)) =0$. 
\item A general member of the linear system $|-K_Y|$ is irreducible and reduced.
\end{enumerate}
\end{prop}

\begin{proof}
Taking the minimal resolution of $Y$, we may assume that $Y$ is smooth. 
Using Proposition \ref{p-cohomology-can-dP}, the same proof of \cite[Theorem 8.3.2.i]{Dol12} works in our setting, so that (1) holds and general members of $|-K_Y|$ are irreducible. 

It is enough to show that a general member of $|-K_Y|$ is reduced. 
Suppose it is not. 
Then there exist $a>1$ such that a general member is of the form $aC \in |-K_Y|$ for some curve $C$. 
In particular, $C$ is a smooth rational curve by Lemma \ref{l-comp-antican-sys}.
Recall that we have the short exact sequence
\[ 0 \rightarrow H^0(Y, \mathcal{O}_Y) \rightarrow  H^0(Y, \mathcal{O}_Y(C)) \rightarrow  H^0(C, \mathcal{O}_C(C)) \rightarrow 0. \]
Since $H^1(Y, \MO_Y)=0$ (Proposition \ref{p-cohomology-can-dP}), 
we have that $h^0(Y, \mathcal{O}_Y(C)) = 1+h^0(C, \mathcal{O}_C(C))$. As $C$ is a smooth rational curve, we conclude by the Riemann--Roch theorem that $h^0(Y, \mathcal{O}_Y(C)) = 2 + C^2$. 

We now consider the induced map 
\begin{eqnarray*}
\theta:H^0(Y, \MO_Y(C) ) &\to& H^0(Y, \MO_Y(aC) ) \simeq H^0(Y, \MO_Y(-K_Y) )\\
\varphi & \mapsto & \varphi^a. 
\end{eqnarray*}
Since a general member of $|-K_Y|$ is of the form $aD$ for some $D \geq 0$, 
$\theta$ is a dominant morphism if we consider $\theta$ as a morphism of affine spaces. 
Therefore, it holds that 
\[ 
h^0(Y, \mathcal{O}_Y(-K_Y)) 
\leq h^0(Y, \mathcal{O}_Y(C)) =2+C^2 =-K_Y\cdot C \leq K_Y^2, \]
which contradicts Proposition \ref{p-cohomology-can-dP}.
\end{proof}

\subsubsection{Anti-canonical systems on geometrically canonical del Pezzo surfaces}

In this section, we study anticanonical systems on geometrically canonical del Pezzo surfaces over an arbitrary field $k$ and we describe their anti-canonical model when the anticanonical degree is small. 

We need the following results on geometrically integral curves of genus one.
\begin{lem} \label{l-gen-genus-one}
Let $k$ be a field. 
Let $C$ be a geometrically integral Gorenstein projective curve over $k$ of arithmetic genus one with $k=H^0(C, \mathcal{O}_C)$. 
Let $L$ be a Cartier divisor on $C$ and let $R(C,L):= \bigoplus_{m \geq 0} H^0(C, mL)$ be the graded $k$-algebra. Then the following hold.
\begin{enumerate}
\item [(i)] If $\deg_k (L) =1$, then $\Bs(L)=\{P\}$ for some $k$-rational point $P$ and $R(C, L)$ is generated 
by $\bigoplus_{1 \leq j \leq 3} H^0( C, jL )$ as a $k$-algebra. 
\item [(ii)] If $\deg_k(L) \geq 2$, then $L$ is globally generated and $R(C, L)$ is generated 
by $H^0( C, L ) \oplus H^0(C, 2L)$ as a $k$-algebra. 
\item [(iii)] If $\deg_k L \geq 3$, then $L$ is very ample and $R(C, L)$ is generated by $H^0( C, L )$ as a $k$-algebra. 
\end{enumerate}
\end{lem}
\begin{proof}
See \cite[Lemma 11.10 and Proposition 11.11]{Tan19}. 
\end{proof}

\begin{prop}\label{p-antican-ring-dP}
Let $k$ be a field. 
Let $X$ be a geometrically canonical weak del Pezzo surface over $k$ 
such that $k=H^0(X, \mathcal{O}_X)$.
Let $R(X, -K_X) = \bigoplus_{m \geq 0} H^0(X, \mathcal{O}_X(-mK_X))$ be the graded $k$-algebra. Then the following hold. 
\begin{enumerate}
\item If $m$ is a positive integer such that $mK_X^2 \geq 2$, 
then $|-mK_X|$ is base point free. 
\item If $K_X^2=1$, then $\Bs(-K_X)=\{P\}$ for some $k$-rational point $P$.
\item If $K_X^2=1$, then $R(X, -K_X)$ is generated by $\bigoplus_{1 \leq j \leq 3} 
H^0(X, -jK_X)$ as a $k$-algebra. 
\item If $K_X^2 =2$, then $R(X, -K_X)$ is generated by $H^0(X, -K_X) \oplus H^0(X, -2K_X)$ 
as a $k$-algebra. 
\item If $K_X^2 \geq 3$, then $R(X, -K_X)$ is generated by $H^0(X, -K_X)$ 
as a $k$-algebra. 
\end{enumerate}
In particular, if $-K_X$ is ample, then $|-6K_X|$ is very ample.
\end{prop}

\begin{proof}
Consider the following condition. 
\begin{enumerate}
\item[(2)'] If $K_X^2=1$, then $\Bs(-K_X)$ is not empty and of dimension zero.  
\end{enumerate}
Since $K_X^2=1$, (2) and (2)' are equivalent. 
Note that to show that (1), (2)', and (3)--(5), 
we may assume that $k$ is algebraically closed. 

From now on, let us prove (1)--(5) under the condition that $k$ is algebraically closed. 
It follows from Proposition \ref{p-gen-memb-antican} 
that a general member $C$ of $|-K_X|$ is a prime divisor. 

Since $C$ is a Cartier divisor and $X$ is Gorenstein, 
then $C$ is a Gorenstein curve. 
By adjunction, $C$ is a Gorenstein curve of arithmetic genus $p_a(C)=1$. 
By Proposition \ref{p-cohomology-can-dP}, we have the following exact sequence for every integer $m$:
\[ 0 \to  H^0(X, -(m-1)K_X) \to H^0(X,-mK_X) \to H^0(C, -mK_X|_C)  \to 0. \]
By the above exact sequence, 
the assertions (1) and (2) follow 
from (3) and (2) of Lemma \ref{l-gen-genus-one}, respectively. 

We prove the assertions (3), (4) and (5).
By the above short exact sequence, it is sufficient to prove the same statement for the $k$-algebra $R(C, \mathcal{O}_C(-K_X))$, which is the content of Lemma \ref{l-gen-genus-one}.
\end{proof}

\begin{thm}\label{t-dP-small-degree}
Let $k$ be a field. 
Let $X$ be a geometrically canonical del Pezzo surface over $k$ such that $H^0(X, \MO_X)=k$.  
Then the following hold. 
\begin{enumerate}
\item If $K_X^2=1$, then $X$ is isomorphic to a weighted hypersurface in $\mathbb{P}_k(1,1,2,3)$ 
of degree six.
\item If $K_X^2=2$, then 
 $X$ is isomorphic to a weighted hypersurface in $\mathbb{P}_k(1,1,1,2)$ 
of degree four. 
\item If $K_X^2=3$, then $X$ is isomorphic to a hypersurface in $\mathbb{P}_k^3$ 
of degree three. 
\item If $K_X^2=4$, then $X$ is isomorphic to a complete intersection of two quadric hypersurfaces in $\mathbb{P}_k^4$.
\end{enumerate}
\end{thm}

\begin{proof}
Using Proposition \ref{p-antican-ring-dP}, the proof is the same as in \cite[Theorem III.3.5]{Kol96}.
\end{proof}

\subsection{Mori fibre spaces to curves}
In this subsection, we summarise properties of regular curves with anti-ample canonical divisor and of Mori fibre space of dimension two over arbitrary fields.

\begin{lem}\label{l-conic}
Let $k$ be a field. 
Let $C$ be a projective Gorenstein integral curve over $k$. 
Then the following are equivalent. 
\begin{enumerate}
\item $\omega_C^{-1}$ is ample. 
\item $H^1(C, \MO_C)=0$. 
\item $C$ is a conic curve of $\mathbb P^2_K$, where $K:=H^0(C, \MO_C)$. 
\item $\deg_k \omega_C = -2 \dim_k (H^0(C, \MO_C))$. 
\end{enumerate}
\end{lem}

\begin{proof}
It follows from \cite[Corollary 2.8]{Tan18a} 
that (1), (2), and (4) are equivalent. 
Clearly, (3) implies (1). 
By \cite[Lemma 10.6]{Kol13}, (1) implies (3). 
\end{proof}

\begin{lem}\label{p-Fano-curve}
Let $k$ be a field and let $C$ be a projective Gorenstein integral curve over $k$ 
such that $k=H^0(C, \MO_C)$ and $\omega_C^{-1}$ is ample. 
Then the following hold. 
\begin{enumerate}
\item If $C$ is geometrically integral over $k$, then $C$ is smooth over $k$. 
\item If the characteristic of $k$ is not two, then $C$ is geometrically reduced over $k$. 
\item If the characteristic of $k$ is not two and $C$ is regular, then $C$ is smooth over $k$.
\end{enumerate}
\end{lem}

\begin{proof}
By Lemma \ref{l-conic}, $C$ is a conic curve in $\mathbb{P}^2_k$. 
Thus, the assertion (1) follows from the fact that an integral conic curve over an algebraically closed field is smooth. 

Let us show (2) and (3). 
Since the characteristic of $k$ is not two and $C$ is a conic curve in $\mathbb P^2_k$, 
we can write 
\[
C=\Proj\,k[x, y, z]/(ax^2+by^2+cz^2)
\]
for some $a, b, c \in k$. 
Since $C$ is an integral scheme, two of $a, b, c$ are not zero. 
Hence, $C$ is reduced. 
Thus (2) holds. 
If $C$ is regular, then each of $a, b, c$ is nonzero, hence $C$ is smooth over $k$. 
\end{proof}

\begin{prop}\label{p-MFS-basic}
Let $k$ be a field. 
Let $\pi:X \to B$ be a $K_X$-Mori fibre space 
from a projective regular $k$-surface $X$ to a projective regular $k$-curve with $k=H^0(B, \MO_B)$. 
Let $b$ be a (not necessarily closed) point. 
Then the following hold. 
\begin{enumerate}
\item The fibre $X_b$ is irreducible.
\item The equation $\kappa(b)=H^0(X_b, \MO_{X_b})$ holds.
\item The fibre $X_b$ is reduced.  
\item The fibre $X_b$ is a conic in $\mathbb P^2_{\kappa(b)}$.  
\item If ${\rm char}\,k \neq 2$, then any fibre of $\pi$ is geometrically reduced. 
\item If ${\rm char}\,k \neq 2$ and $k$ is separably closed, 
then $\pi$ is a smooth morphism. 
\end{enumerate}
\end{prop}

\begin{proof}
If $X_b$ is not irreducible, it contradicts the hypothesis $\rho(X/B)=1$. 
Thus (1) holds. 

Let us show (2). 
Since $\pi$ is flat, the integer
\[
\chi:= \dim_{\kappa(b)}H^0(X_b, \MO_{X_b}) - \dim_{\kappa(b)}H^1(X_b, \MO_{X_b}) \in \Z
\]
is independent of $b\in B$. 
Since $H^1(X_b, \MO_{X_b})=0$ for any $b \in B$, it suffices to show that 
$\dim_{\kappa(b)}H^0(X_b, \MO_{X_b})=1$ for some $b \in B$. 
This holds for the case when $b$ is the generic point of $B$. 
Hence, (2) holds. 

Let us prove (3). 
It is clear that the generic fibre is reduced. 
We may assume that $b \in B$ is a closed point. 
Assume that $X_b$ is not reduced. 
By (1), we have $X_b=mC$ for some prime divisor $C$ and $m \in \Z_{\geq 2}$. 
Since $-K_X \cdot_{\kappa(b)} X_b=2$, we have that $m=2$. 
Then we obtain an exact sequence: 
\[
0 \to \MO_X(-C)|_C \to \MO_{X_b} \to \MO_C \to 0.
\]
Since $C^2=0$ and $\omega_{C}^{-1}$ is ample, we have that $\MO_X({-C})|_C \simeq \MO_C$. 
Since $H^1(C, \MO_C)=0$, we get an exact sequence: 
\[
0 \to H^0(C, \MO_C) \to H^0(X_b, \MO_{X_b}) \to H^0(C, \MO_C) \to 0. 
\] 
Then we obtain $\dim_{\kappa(b)} H^0(X_b, \MO_{X_b}) \geq 2$, 
which contradicts (2). 
Hence (3) holds. 

We now show (4). By \cite[Corollary 2.9]{Tan18a},
$\deg_{\kappa(b)} \omega_{X_b}=(K_X+X_b) \cdot_{\kappa(b)} X_b<0$. Hence (4) follows from (2) and Lemma \ref{l-conic}. 

The assertions (5) and (6) follow from Proposition \ref{p-Fano-curve}. 
\end{proof}

\subsection{Twisted forms of canonical singularities}

The aim of this subsection is to prove Proposition \ref{p-insep-bdd-rat-pts}. 
The main idea is to bound the purely inseparable degree of regular non smooth points on geometrically normal surfaces according to the type of singularities.
For this, the notion of Jacobian number plays a crucial role.

\begin{dfn} \label{d-number}
Let $k$ be a field of characteristic $p>0$. 
Let $R$ be an equi-dimensional $k$-algebra essentially of finite type over $k$. 
Let $J_{R/k}$ be its Jacobian ideal of $R$ over $k$ (cf. \cite[Definition 4.4.1 and Proposition 4.4.4]{HS06}). 
We define the {\em Jacobian number} of $R/k$ as $\nu(R):=\nu(R/k) := \dim_k (R/J_{R/k})$. 
Note that $\nu(R/k) <\infty$ if $R/J_{R/k}$ is an artinian ring 
and its residue fields are finite extensions of $k$. 
\end{dfn}

\begin{rem}\label{r-base-ch}
Let $k \subset k'$ be a field extension of characteristic $p>0$ 
and let $R$ be an equi-dimensional $k$-algebra essentially of finite type over $k$. 
Then the following hold. 
\begin{enumerate}
\item 
By \cite[Definition 4.4.1]{HS06}, we get 
\[
J_{R/k} \cdot (R \otimes_k k') = J_{R \otimes_k k'/k'}.
\]
In particular, if $R/J_{R/k}$ is an artinian ring and its residue fields are finite extensions of $k$, 
then we have $\nu(R/k)=\nu(R \otimes_k k'/k')$.  
\item 
Assume that $k$ is a perfect field. 
By \cite[Definition 4.4.9]{HS06}, 
$\Spec\,(R/J_{R/k})$ set-theoretically coincides with 
the non-regular locus of $\Spec\,R$.  
\item 
Assume that $R$ is of finite type over $k$. 
Then (1) and (2) imply that $\Spec\,(R/J_{R/k})$ set-theoretically coincides with 
the the non-smooth locus of $\Spec\,R \to \Spec\,k$.
\end{enumerate}
\end{rem}

\begin{rem}\label{r-nu-2dim-gn}
In our application, 
$R$ will be assumed to be a local ring $\MO_{X, x}$ at a closed point $x$ 
of a geometrically normal surface $X$ over $k$. 
In this case, (3) of Remark \ref{r-base-ch} implies that 
$R/J_{R/k}$ is an artinian local ring 
whose residue field is a finite extension of $k$. 
Hence, $\nu(R/k) = \dim_k (R/J_{R/k})$ is well-defined as in Definition \ref{d-number}. 
\end{rem}

To treat local situations, 
let us recall the notion of essentially \'etale ring homomorphisms.  
For its fundamental properties, we refer to \cite[Subsection 2.8]{Fu15}. 

\begin{dfn}\label{d-ess-et}
Let $f\colon R \to S$ be a local homomorphism of local rings. 
We say that $f$ is {\em essentially \'etale} 
if there exists an \'etale $R$-algebra $\overline S$ and a prime ideal $\mathfrak p$ 
of $\overline S$ such that 
$\mathfrak p$ lies over the maximal ideal of $R$ and $S$ is $R$-isomorphic to $\overline S_{\mathfrak p}$. 
\end{dfn}

\begin{lem}\label{l-nu-et}
Let $k$ be a field.
Let $f \colon R \to S$ be an essentially \'etale local $k$-algebra homomorphism 
of local rings which are essentially of finite type over $k$. 
Let $\m_R$ and $\m_S$ be the maximal ideals of $R$ and $S$, respectively. 
Set $\kappa(R):=R/\m_R$ and $\kappa(S):=S/\m_S$. 
Then the following hold. 
\begin{enumerate}
\item 
If $M$ is an $R$-module of finite length whose support is contained the maximal ideal $\m_R$, 
then the equation 
\[
\dim_k (M \otimes_R S) = [\kappa(S):\kappa(R)] \dim_k M
\]
 holds. 
\item 
Suppose that $R$ is an integral domain, 
$R/J_{R/k}$ is an artinian ring, and $\kappa(R)$ is a finite extension of $k$.
Then the equation 
\[
\nu(S/k)=[\kappa(R):\kappa(S)]\nu(R/k)
\]
holds. 
\end{enumerate}
\end{lem}

\begin{proof}
Let us show (1). 
Since $M$ is a finitely generated $R$-module,
there exists a sequence of $R$-submodules $M=:M_0 \supset M_1 \supset \cdots \supset M_n =0$ 
such that $M_i/M_{i+1} \simeq R/ \p$ for some prime ideal $\p$ by \cite[Theorem 6.4]{Mat89} . 
Since the support of $M$ is $\m_R$, we have $\p=\m_R$. 
As $R \to S$ is flat, the problem is reduced to the case when $M=R/\m_R=\kappa(R)$. 
In this case, we have 
\[
\kappa(R) \otimes_R S = (R/\m_R) \otimes_R S \simeq S/\m_R S =S/\m_S = \kappa(S), 
\]
where the equality $S/m_R S = S/m_S$ follows from the assumption that $f$ is a localisation of 
an unramified homomorphism. 
Hence, (1) holds. 

Let us show (2). 
Set $n:=\dim R$. 
We use the description of the Jacobian of $R$ via Fitting ideals 
(cf. \cite[Discussion 4.4.7]{HS06}): 
$J_{R/k}= \text{Fit}_n (\Omega_{R/k}^1)$ 
and 
$J_{S/k}= \text{Fit}_n (\Omega_{S/k}^1)$. 
We have 
\[
J_{S/k}= \text{Fit}_n (\Omega_{S/k}^1)= \text{Fit}_n (\Omega^1_{R/k} \otimes_R S)=\text{Fit}_n (\Omega^1_{R/k})S 
=J_{R/k}S,
\]
where the third equality follows from (3) of \cite[\href{https://stacks.math.columbia.edu/tag/07ZA}{Tag 07ZA}]{StackProject}. 
As $f: R \to S$ is flat, we obtain $S/J_{S/k} \simeq (R/J_{R/k}) \otimes_R S$. 
By (1) and Definition \ref{d-number}, the assertion (2) holds. 
\end{proof}

\begin{ex}\label{ex-A_p^n}
Let $k$ be a field of characteristic $p>0$. Let $X=\Spec\,R$ be a surface over $k$ such that
\begin{enumerate}
\renewcommand{\labelenumi}{(\roman{enumi})}
\item $X \times_k {\overline{k}}=\Spec\,(R\otimes_k \overline k)$ is a normal surface, 
\item  $X \times_k {\overline{k}}$ has a unique singular point $x$, 
and $x$ is a canonical singularity of type $A_{p^n-1}$. 
\end{enumerate}
We prove that $\nu(R/k) = p^n$.
By Remark \ref{r-base-ch}, we have $\nu(R/k)=\nu(R \otimes_k \overline{k}/\overline{k}).$ 
In order to compute $\nu(R \otimes_k \overline{k}/\overline{k})$, it is sufficient to localise at the singular point by \cite[Corollary 4.4.5]{HS06}.
Thus we can suppose that $k$ is algebraically closed and $R$ is a local $k$-algebra.

By \cite[pages 16-17]{Art77} (cf. (12) of Subsection \ref{ss-notation}), 
the henselisation $R^h$ of $R$ is isomorphic to
\[ k\{x,y,z\}/(z^{p^n}+xy). \] 
In particular there exist essentially \'etale local $k$-algebra homomorphisms 
$R \to S$ and $k[x,y,z]/(z^{p^n}-xy) \to S.$
A direct computation shows $\nu(k[x,y,z]/(z^{p^n}-xy))=p^n$.
Thus by Lemma \ref{l-nu-et}, we have
\[ \nu(R) = \nu (S) = \nu(k[x,y,z]/(z^{p^n}-xy))=p^n.  \]
\end{ex}

The following is a generalisation of \cite[Lemma 14.2]{FS18}.

\begin{lem}\label{l-bound-points}
Let $k$ be a field of characteristic $p>0$. 
Let $X=\Spec\, R$, where $R$ is an equi-dimensional local $k$-algebra of essentially finite type over $k$. 
Let $x$ be the closed point of $X$. 
Suppose that 
$R/J_{R/k}$ is a local artinian ring and 
its residue field $\kappa(x)$ is a finite extension of $k$. 
Then  $[ \kappa(x) : k]$ is a divisor of $\nu(R/k)$.
\end{lem}

\begin{proof} 
Let $R/J_{R/k}=:M_0 \supset M_1 \supset \cdots \supset M_n =0$ be a composition sequence of 
$R/J_{R/k}$-submodules (cf. \cite[Theorem 6.4]{Mat89}). 
Since $R/J_{R/k}$ is an artinian local ring, it holds that $M_i/M_{i+1} \simeq \kappa(x)$ for any $i$. 
We have 
{\small 
\[
\nu(R/k)=
\dim_k (R/J_{R/k}) = \sum_{i=0}^{n-1} \dim_k (M_i/M_{i+1})
=n \dim_k \kappa(x) = n[\kappa(x):k]. 
\]
}
We thus conclude that $[\kappa(x):k]$ is a divisor of $\nu(R/k)$. 
\end{proof}

\begin{lem} \label{l-deg-ext-sing}
Let $X$ be a regular variety over a separably closed field $k$. 
Suppose that $X_{\overline{k}} = X \times_k \overline k$ is a normal variety with a unique singular point $y$. 
Let $x$ be the image of $y$ by the induced morphism $X_{\overline k} \to X$. 
Then the following hold. 
\begin{enumerate}
\item $[\kappa(x) : k]$ is a divisor of $\nu(\MO_{X,x})$. 
\item $X \times_k \kappa(x)$ is not regular.
\end{enumerate}
\end{lem}

\begin{proof}
Since $k$ is separably closed, 
the induced morphism $X_{\overline{k}} \to X$ is a universal homeomorphism. 
Note that the local ring $\mathcal{O}_{X,x}$ is not geometrically regular over $k$.
Applying Lemma \ref{l-bound-points} to the local ring $\mathcal{O}_{X,x}$, 
we deduce that  
$[\kappa(x):k]$ is a divisor of $\nu(\MO_{X,x})$. 
Thus (1) holds. 
Consider the base change $\pi \colon X \times_k \kappa(x) \rightarrow X$. 
Let $x'$ be the point on $X \times_k \kappa(x)$ lying over $x$. 
Note that $x'$ is a $\kappa(x)$-rational point of $X \times_k \kappa(x)$ whose base change by $(-) \times_{\kappa(x)} \overline k$ is not regular. 
By \cite[Corollary 2.6]{FS18}, 
we conclude that $X \times_k \kappa(x)$ is not regular at $x'$.
\end{proof}

We now explain how the previous results can be used to construct closed points with purely inseparable residue field on a regular surface. This will be used in Section \ref{s-pi-pts} to find purely inseparable points on regular del Pezzo surfaces.

\begin{prop} \label{p-insep-bdd-rat-pts}
Let $X$ be a regular surface over $k$. 
Suppose that $X_{\overline k}=X \times_k \overline k$ is a normal surface over $\overline k$ 
with a unique singular point $y$. 
Assume that $y$ is a canonical singularity of type $A_{p^n-1}$.
Let $z$ be the image of $y$ by the induced morphism $X_{\overline k} \to X_{k^{1/p^n}} =X \times_k k^{1/p^n}$. 
Then $z$ is a $k^{1/p^n}$-rational point on $X_{k^{1/p^n}}$. 
\end{prop}

\begin{proof}
Set $R:=\MO_{X, x}$, where $x$ is the unique closed point along which $X$ is not smooth.
Let $k^{\text{sep}}$ be the separable closure of $k$. 
For $R_{k^{\text{sep}}}:=R \otimes_k k^{\text{sep}}$, 
it follows from Example \ref{ex-A_p^n} that $\nu(R_{k^{\text{sep}}})=p^n$. 
Lemma \ref{l-deg-ext-sing} implies 
that $k^{\sep} \subset \kappa(z)$ is purely inseparable and $[\kappa(z):k]$ is a divisor of $p^n$. 
In particular, $\kappa(z) \subset (k^{\text{sep}})^{1/p^n}$.

Consider the Galois extension $k^{1/p^n} \subset (k^{\text{sep}})^{1/p^n}$ and denote by  $G$ its Galois group. 
For $X_{(k^{\text{sep}})^{1/p^n}}:=X \times_k (k^{\text{sep}})^{1/p^n}$, 
$G$ acts on the set $X_{(k^{\text{sep}})^{1/p^n}}((k^{\text{sep}})^{1/p^n})$.  
The unique singular $(k^{\text{sep}})^{1/p^n}$-rational point on   $X_{(k^{\text{sep}})^{1/p^n}}$ 
is fixed under the $G$-action.
Thus it descends to a $k^{1/p^n}$-rational point on $X_{k^{1/p^n}}$.
\end{proof}

\section{Behaviour of del Pezzo surfaces under base changes}

In this section, we study the behaviour of canonical del Pezzo surfaces over an imperfect field $k$ under the base changes to the algebraic closure $\overline{k}$. 

\subsection{Classification of base changes of del Pezzo surfaces}\label{s-classify}

In this subsection, we give classification of base changes of 
del Pezzo surfaces with canonical singularities over imperfect fields 
(Theorem \ref{t-classify-bc}). 
To this end, we need two auxiliary lemmas: 
Lemma \ref{l-Reid} and Lemma \ref{l-rationality}. 
The former one classify $\Q$-factorial surfaces over algebraically closed fields 
whose anti-canonical bundles are sufficiently positive. 
Its proof is based on a simple but smart idea by Reid 
(cf. the proof of \cite[Theorem 1.1]{Rei94}). 
The latter one, i.e. Lemma~\ref{l-rationality}, 
gives a rationality criterion 
for the base changes of log del Pezzo surfaces.

\begin{lem}\label{l-Reid}
Let $k$ be an algebraically closed field. 
Let $Y$ be a projective normal $\Q$-factorial surface over $k$ 
such that $-K_Y \equiv A+D$ 
for an ample Cartier divisor $A$ and a pseudo-effective $\Q$-divisor $D$. 
Let $\mu:Z \to Y$ be the minimal resolution of $Y$. 
Then one of the following assertions holds. 
\begin{enumerate}
\item $D\equiv 0$ and $Y$ has at worst canonical singularities. 
\item $Z$ is 
isomorphic to a $\mathbb P^1$-bundle over a smooth projective curve. 
\item $Z \simeq \mathbb P^2$. 
\end{enumerate} 
\end{lem}

\begin{proof}
Assuming that (1) does not hold, let us prove that either (2) or (3) holds. 
We have 
\[
K_Z+E=\mu^*K_Y
\]
for some effective $\mu$-exceptional $\Q$-divisor $E$ on $Z$. 
In particular, it holds that 
\[
K_Z+E+\mu^*(D)=\mu^*(K_Y+D) \equiv -\mu^*A. 
\]
Since (1) does not hold, we have that $D \not\equiv 0$ or $E \neq 0$. 
Then we get 
\[
K_Z+\mu^*A \equiv -E-\mu^*(D) \not\equiv 0, 
\]
hence $K_Z+\mu^*A$ is not nef. 
By the cone theorem for a smooth projective surface \cite[Theorem 1.24]{KM98}, 
there is a curve $C$ that spans  
a $(K_Z+\mu^*A)$-negative extremal ray $R$ of $\overline{{\rm NE}}(Z)$. 
Note that $C$ is not a $(-1)$-curve. 
Indeed, otherwise $\mu(C)$ is a curve and we obtain $\mu^*A \cdot C>0$, 
which induces a contradiction: 
\[
(K_Z+\mu^*A) \cdot C\geq -1+1=0. 
\] 
It follows from the classification of the $K_Z$-negative extremal rays 
\cite[Theorem 1.28]{KM98} that 
either $Z \simeq \mathbb P^2$ or 
$Z$ is a $\mathbb P^1$-bundle over a smooth projective curve. 
In any case, one of (2) and (3) holds.  
\end{proof}

\begin{lem}\label{l-rationality}
Let $(X, \Delta)$ be a projective two-dimensional klt pair over a field of characteristic $p>0$ such that $-(K_X+\Delta)$ is nef and big. 
Assume that $k=H^0(X, \MO_X)$. 
Then $(X \times_k \overline{k})_{\red}$ is a rational surface. 
\end{lem}

\begin{proof}
See \cite[Proposition 2.20]{NT}. 
\end{proof}

We now give a classification of the base changes of del Pezzo surfaces with canonical singularities.

\begin{thm}\label{t-classify-bc}
Let $k$ be a field of characteristic $p>0$. 
Let $X$ be a canonical del Pezzo surface over $k$ with $k=H^0(X, \MO_X)$. 
Then the normalisation $Y$ of $(X \times_k \overline{k})_{\red}$ satisfies one of the following properties. 
\begin{enumerate}
\item $X$ is geometrically canonical over $k$. 
In particular, $Y \simeq X \times_k \overline k$  and $-K_Y$ is ample. 
\item $X$ is not geometrically normal over $k$ and 
$Y$ is isomorphic to a Hirzebruch surface, i.e. a $\mathbb P^1$-bundle over $\mathbb P^1$. 
\item $X$ is not geometrically normal over $k$ and $Y$ is isomorphic to a weighted projective surface $\mathbb P(1, 1, m)$ 
for some positive integer $m$. 
\end{enumerate} 
\end{thm}

\begin{proof}
Replacing $k$ by its separable closure, 
we may assume that $k$ is separably closed. 
Let $f:Y \to X$ be the induced morphism and let $\mu:Z \to Y$ be the minimal resolution of $Y$. 
By \cite[Theorem 4.2]{Tan18b}, 
there is an effective $\Z$-divisor $D$ on $Y$ such that 
\begin{itemize}
\item $K_Y+D=f^*K_X$, and 
\item if $X \times_k \overline k$ is not normal, then $D \neq 0$. 
\end{itemize}
Since $-K_X$ is an ample Cartier divisor, so is $-f^*K_X$. 
Moreover, 
it follows from \cite[Lemma 2.2 and Lemma 2.5]{Tan18b} 
that $Y$ is $\Q$-factorial. 
Hence, we may apply Lemma~\ref{l-Reid} to $-K_Y= -f^*K_X+D$.

By Lemma~\ref{l-rationality}, $Y$ is a rational surface. 
Thus, if (2) or (3) of Lemma~\ref{l-Reid} holds, 
then one of (1)--(3) of Theorem \ref{t-classify-bc} holds, as desired. 
Therefore, let us treat the case when (1) of Lemma~\ref{l-Reid} holds. 
Then it holds that $D=0$ and $Y$ has at worst canonical singularities. 
In this case, we have that $Y= X \times_k \overline k$ and $X$ is geometrically canonical. 
Hence, (1) of Theorem \ref{t-classify-bc} holds, as desired. 
\end{proof}

\subsection{Bounds on Frobenius length of geometric non-normality}

In this subsection, we give an upper bound for 
the Frobenius length of geometric non-normality for canonical del Pezzo surfaces 
(Proposition \ref{p-p2-bound}). 
We start by recalling its definition (Definition \ref{d-lF}) and 
fundamental properties (Remark \ref{r-lF}).

\begin{dfn}\label{d-lF}
Let $k$ be a field of characteristic $p>0$. 
Let $X$ be a proper normal variety over $k$ such that $k=H^0(X, \MO_X)$. 
The {\em Frobenius length of geometric non-normality}  $\ell_F(X/k)$  of $X/k$ 
is defined by 
{\small 
\[
\ell_F(X/k):=\min\{\ell \in \Z_{\geq 0}\,|\, 
(X \times_k k^{1/p^{\ell}})_{\red}^N \text{ is geometrically normal over }k^{1/p^{\ell}}\}.
\]
}
\end{dfn}

\begin{rem}\label{r-lF} 
Let $k$ and $X$ be as in Definition \ref{d-lF}. 
Set $\ell:=\ell_F(X/k)$. Let $(k', Y)$ be one of $(k^{1/p^{\infty}}, (X \times_k k^{1/p^{\infty}})_{\red}^N)$ 
and  $(\overline k, (X \times_k \overline k)_{\red}^N)$. 
We summarise some results from \cite[Section 5]{Tan19}. 

\begin{enumerate}
\item 
The existence of the right hand side of Definition \ref{d-lF} 
is assured by \cite[Remark 5.2]{Tan19}. 
\item 
If $X$ is not geometrically normal, 
then $\ell$ is a positive integer \cite[Remark 5.3]{Tan19} and 
there exist nonzero effective Weil divisors $D_1, ..., D_{\ell}$ 
such that 
\[
K_Y+(p-1)\sum_{i=1}^{\ell} D_i \sim f^*K_X,
\]
where $f:Y \to X$ denotes the induced morphism \cite[Proposition 5.11]{Tan19}. 
\item 
The $\ell$-th iterated absolute Frobenius morphism $F^{\ell}_{X\times_k k'}$ factors through 
the induced morphism $Y \to X \times_k k'$ \cite[Proposition 5.4 and Theorem 5.9]{Tan19}: 
\[
F^{\ell}_{X\times_k k'}:X \times_k k' \to Y \to X \times_k k'. 
\]
\end{enumerate}
\end{rem}

\begin{prop}\label{p-p2-bound}
Let $k$ be a field of characteristic $p>0$. 
Let $X$ be a canonical del Pezzo surface over $k$ with $k=H^0(X, \MO_X)$. 
Let $Y$ be the normalisation of $(X \times_k \overline k)_{\red}$ and let $f:Y \to X$ be the induced morphism. 
Assume that the linear equivalence 
\[
K_Y+\sum_{i=1}^r C_i \sim f^*K_X
\]
holds for some prime divisors $C_1, ..., C_r$ (not necessarily $C_i \neq C_j$ for $i \neq j$). 
Then it holds that $r \leq 2$. 
\end{prop}

\begin{proof}
Set $C:=\sum_{i=1}^r C_i$. We have $K_Y+C \sim f^*K_X$. 
If $C=0$, then there is nothing to show. 
Hence, we may assume that $C \neq 0$. 
In particular, $X$ is not geometrically normal. 
In this case, it follows from Theorem \ref{t-classify-bc} 
that $Y$ is isomorphic to either a Hirzebruch surface or $\mathbb P(1, 1, m)$ for some $m>0$.

We first treat the case when $Y \simeq \mathbb P(1, 1, m)$. 
If $m=1$, then the assertion is obvious. 
Hence, we may assume that $m \geq 2$. 
In this case, for the minimal resolution $g:Z \to Y$, 
we have that 
\[
K_Z+\frac{m-2}{m}\Gamma=g^*K_Y
\]
where $\Gamma$ is the negative section of the fibration $ Z \rightarrow \mathbb{P}^1$ 
such that $\Gamma^{2}=-m$. Note that $m$ is the $\Q$-factorial index of $Y$, i.e. $m D$ is Cartier for any $\Z$-divisor $D$ on $Y$. 
We have that 
$$-K_Z=\frac{m-2}{m}\Gamma-g^*K_Y\equiv \frac{m-2}{m}\Gamma+g^*C-g^*f^*K_X$$
Consider the intersection number with 
a fibre $F_Z$ of $Z \to \mathbb P^1$: 
$$2=\left(\frac{m-2}{m}\Gamma+g^*C-g^*f^*K_X\right) \cdot F_Z 
\geq \frac{m-2}{m}+C \cdot g_*(F_Z) +1.$$
Thus we obtain 
$$2 \geq C \cdot (mg_*(F_Z)) \geq r,$$
where the last inequality holds since $mg_*(F_Z)$ is an ample Cartier divisor. 
Therefore, we obtain $r \leq 2$, as desired. 

It is enough to treat the case when $Y$ is a Hirzebruch surface. 
For a fibre $F$ of $\pi: Y \to \mathbb P^1$, we have that 
\[
-2+C \cdot F= (K_Y+C) \cdot F = f^*K_X \cdot F \leq -1,
\]
hence $C \cdot F \leq 1$. 
There are two possibilities: $C \cdot F=1$ or $C \cdot F=0$. 

Assume that $C \cdot F=1$. 
Then there is a section $\Gamma$ of $\pi$ and a $\pi$-vertical $\Z$-divisor $C'$ such that $C=\Gamma+C'$. 
Consider the intersection number with $\Gamma$: 
\[
-2+\Gamma \cdot C'=(K_Y+\Gamma+C') \cdot \Gamma =(K_Y+C) \cdot \Gamma = f^*K_X \cdot \Gamma \leq -1. 
\]
Therefore, we have $\Gamma \cdot C' \leq 1$. 
This implies that either $C'=0$ or $C'$ is a prime divisor. 
In any case, we get $r \leq 2$, as desired.

We may assume that $C \cdot F=0$, i.e. $C$ is a $\pi$-vertical divisor. 
Let $\Gamma$ be a section of $\pi$ such that $\Gamma^2 \leq 0$. 
We have that 
\[
-2 + C \cdot \Gamma=(K_Y+\Gamma+C) \cdot \Gamma  \leq (K_Y+C) \cdot \Gamma =f^*K_X \cdot \Gamma \leq -1. 
\]
Hence, we obtain $C \cdot \Gamma \leq 1$, which implies $r \leq 1$. 
\end{proof}

\begin{thm}\label{t-p2-bound}
Let $k$ be a field of characteristic $p>0$. 
Let $X$ be a canonical del Pezzo surface over $k$ 
such that $k=H^0(X, \MO_X)$. 
Let $Y$ be the normalisation of $(X \times_k \overline k)_{\red}$ and let 
\[
\mu: Y \to X \times_k \overline k
\]
be the induced morphism. 
\begin{enumerate}
\item If $p \geq 5$, then $X$ is geometrically canonical, 
i.e. $\mu$ is an isomorphism and $Y$ has at worst canonical singularities. 
\item If $p=3$, then $\ell_{F}(X/k) \leq 1$ and the absolute Frobenius morphism $F_{X \times_k \overline k}$ 
of $X \times_k \overline k$ factors through $\mu$: 
\[
F_{X \times_k \overline k}:X \times_k \overline k\to Y \xrightarrow{\mu} X \times_k \overline k.
\]
\item 
If $p=2$, then $\ell_{F}(X/k) \leq 2$ and the second iterated absolute Frobenius morphism $F^2_{X \times_k \overline k}$ 
of $X \times_k \overline k$ factors through $\mu$: 
\[
F^2_{X \times_k \overline k}:X \times_k \overline k\to Y \xrightarrow{\mu} X \times_k \overline k.
\]
\end{enumerate}
\end{thm}

\begin{proof}
The assertion follows from Remark \ref{r-lF} and Proposition \ref{p-p2-bound}. 
\end{proof}

\section{Numerically trivial line bundles on log del Pezzo surfaces}\label{s-nume-triv}

The purpose of this section is to give an explicit upper bound on the torsion index of numerically trivial line bundles on log del Pezzo surfaces over imperfect fields (Theorem \ref{t-klt-bdd-torsion}). 
To achieve this result, we use the minimal model program to reduce the problem 
to the case when our log del Pezzo surface admits a Mori fibre space structure $\pi:X \to B$. 
The cases $\dim B=0$ and $\dim B=1$ 
will be settled in Theorem \ref{t-cano-bdd-torsion} and 
Proposition \ref{p-ess-klt-bdd-torsion}, respectively.

\subsection{Canonical case}

In this subsection, we study numerically trivial Cartier divisor on del Pezzo surfaces with canonical singularities.
\begin{thm}\label{t-cano-bdd-torsion}
Let $k$ be a field of characteristic $p>0$. 
Let $X$ be a canonical weak del Pezzo surface over $k$ such that $k=H^0(X, \MO_X)$. 
Let $L$ be a numerically trivial Cartier divisor on $X$. 
Then the following hold. 
\begin{enumerate}
\item If $p \geq 5$, then $L \sim 0$. 
\item If $p=3$, then $3L \sim 0$. 
\item If $p=2$, then $4L \sim 0$. 
\end{enumerate}
\end{thm}

\begin{proof}
We first reduce the problem to the case when $-K_X$ is ample. 
It follows from \cite[Theorem 4.2]{Tan18a} that 
$-K_X$ is semi-ample. 
As $-K_X$ is also big, 
$|-mK_X|$ induces a birational morphism 
$f:X \to Y$ to a projective normal surface $Y$. 
Then it holds that $K_Y$ is $\Q$-Cartier and $K_X=f^*K_Y$. 
In particular, $Y$ has at worst canonical singularities. 
Then \cite[Theorem 4.4]{Tan18a} enables us to find a numerically trivial Cartier divisor $L_Y$ on $Y$ such that $f^*L_Y \sim L$.  
Hence the problem is reduced to the case when $-K_X$ is ample.

We only treat the case when $p=2$, as the other cases are easier. 
By Theorem \ref{t-p2-bound}, the second iterated absolute Frobenius morphism 
\[
F^2_{X \times_k \overline k}:X \times_k \overline k \to X \times_k \overline k
\]
factors through the normalisation $(X \times_k \overline k)_{\red}^N$ of $(X \times_k \overline k)_{\red}$: 
\[
F^2_{X \times_k \overline k}:
X \times_k \overline k \to (X \times_k \overline k)_{\red}^N \xrightarrow{\mu}
 X \times_k \overline k, 
\]
where $\mu$ denotes the induced morphism. 
Set $\mathcal L:=\MO_X(L)$ and 
let $\mathcal L_{\overline k}$ be the pullback of $\mathcal L$ to $X \times_k \overline k$. 
Since $(X \times_k \overline k)_{\red}^N$ is a normal rational surface by Lemma \ref{l-rationality}, 
any numerically trivial invertible sheaf is trivial: 
$\mu^*\mathcal L_{\overline k} \simeq \MO_{(X \times_k \overline k)_{\red}^N}$. 
As $F^2_{X \times_k \overline k}$ factors through $\mu$, we have that 
\[
\mathcal L_{\overline k}^4 = 
(F^2_{X \times_k \overline k})^*\mathcal L_{\overline k} \simeq \MO_{X \times_k \overline k}. 
\]
Then it holds that 
\[
H^0(X, \mathcal L^4) \otimes_k \overline k \simeq 
H^0(X \times_k \overline k, \mathcal L^4_{\overline k}) \simeq 
H^0(X \times_k \overline k, \MO_{X \times_k \overline k}) \neq 0. 
\]
Hence we obtain $H^0(X, \mathcal L^4) \neq 0$, i.e. $4L \sim 0$. 
\end{proof}

\subsection{Essential step for the log case}

In this subsection, 
we study the torsion index of numerically trivial line bundles on log del Pezzo surfaces admitting the following special Mori fibre space structure onto a curve.

\begin{nota}\label{n-ess-klt-case}
We use the following notation. 
\begin{enumerate}
\item $k$ is a field of characteristic $p>0$. 
\item $X$ is a regular $k$-surface of del Pezzo type such that 
$k=H^0(X, \MO_X)$ and $\rho(X)=2$. 
\item $B$ is a regular projective curve over $k$ such that $k=H^0(B, \MO_B)$.
\item $\pi:X \to B$ is a $K_X$-Mori fibre space. 
\item 
Let $R = \R_{\geq 0}[\Gamma]$ be the extremal ray which does not correspond to $\pi$, 
where $\Gamma$ denotes a  curve on $X$. 
Note that $\pi(\Gamma)=B$.  
Set $d_{\Gamma}:=\dim_k H^0(\Gamma, \MO_{\Gamma}) \in \Z_{>0}$ and $m_{\Gamma}:=[K(\Gamma):K(B)] \in \Z_{>0}$. 
We denote by $\pi_{\Gamma}:\Gamma \to B$ the induced morphism.   
\item Assume that $K_X \cdot \Gamma >0$. 
\end{enumerate}
\end{nota}

\begin{lem}\label{l-ess-klt-case}
We use Notation \ref{n-ess-klt-case}. 
Then the following hold. 
\begin{enumerate}
\item[(7)] $\Gamma^2 \leq 0$. 
\item[(8)] There exists a rational number $\alpha$ such that $0 \leq \alpha <1$ 
and  $(X, \alpha \Gamma)$ is a log del Pezzo pair.
\end{enumerate}
\end{lem}

\begin{proof}
The assertion (7) follows from Lemma \ref{l-ext-ray} below.  
Let us prove (8). 
By Notation \ref{n-ess-klt-case}(2), there is an effective $\Q$-divisor $\Delta$ such that $(X, \Delta)$ is a log del Pezzo pair. 
We write $\Delta= \alpha \Gamma + \Delta'$ 
for some rational number $0 \leq \alpha<1$ and an effective $\Q$-divisor $\Delta'$ 
with $\Gamma \not\subset \text{Supp}(\Delta')$. 
Since $\overline{\NE}(X)$ is generated by $\Gamma$ and a fibre $F$ of the morphism $\pi:X \to B$, 
we conclude that any prime divisor $C$ such that $C \neq \Gamma$ is nef. 
In particular, $\Delta'$ is nef. 
Hence, $(X, \alpha \Gamma)$ is a log del Pezzo pair. 
Thus, (8) holds. 
\end{proof}

\begin{lem}\label{l-ext-ray}
Let $k$ be a field. 
Let $X$ be a projective $\Q$-factorial normal surface over $k$ 
Let $R=\R_{\geq 0}[\Gamma]$ is an 
extremal ray of $\overline{\NE}(X)$, 
where $\Gamma$ is a curve on $X$. 
If $\Gamma^2 > 0$, then $\rho(X)=1$. 
\end{lem}

\begin{proof}
We may apply the same argument as in 
\cite[Theorem 3.21, Proof of the case where $C^2>0$ in page 20]{Tan14}. 
\end{proof}

The first step is to prove that $m_{\Gamma} \leq 5$ (Proposition \ref{p-cov-deg-bound}). 
To this end, we find an upper bound and a lower bound for $\alpha$ 
(Lemma \ref{l-alpha-upper}, Lemma \ref{l-alpha-lower}).

\begin{lem}\label{l-alpha-upper}
We use Notation \ref{n-ess-klt-case}. 
Take a closed point $b$ of $B$ and set $F_b:=\pi^*(b)$.
Let $\kappa(b)$ be the residue field at $b$ and set $d(b):=[\kappa(b):k]$. 
Then the following hold. 
\begin{enumerate}
\item 
$K_X \cdot_k F_b = -2d(b)$. 
\item 
$\Gamma \cdot_k F_b = m_{\Gamma}d(b)$. 
\item 
If $\alpha$ is a rational number such that $-(K_X+\alpha \Gamma)$ is ample, 
then $\alpha m_{\Gamma}<2$. 
\end{enumerate}
\end{lem}

\begin{proof}
Let us show (1). 
We have that 
\[
\deg_k \omega_{F_b} = (K_X+F_b)\cdot_k F_b= K_X \cdot_k F_b<0.
\]
Hence, Lemma \ref{l-conic} implies that 
\[
K_X \cdot_k F_b=\deg_k \omega_{F_b}=-2d(b). 
\]
Thus (1) holds. 
Clearly, (2) holds. 

Let us show (3). Since $-(K_X+\alpha \Gamma)$ is ample, (1) and (2) imply that 

\[
0> (K_X+\alpha \Gamma) \cdot_k F_b=-2d(b)+\alpha m_{\Gamma}d(b).  
\]
Thus (3) holds. 
\end{proof}

\begin{lem}\label{l-alpha-lower}
We use Notation \ref{n-ess-klt-case}. 
Then the following hold. 
\begin{enumerate}
\item 
$(K_X+\Gamma) \cdot_k \Gamma = -2d_{\Gamma}<0$. 
\item 
For a rational number $\beta$ with $0 \leq \beta \leq 1$, it holds that 
\[
(K_X+\beta \Gamma) \cdot_k \Gamma \geq d_{\Gamma}(1-3\beta). 
\]
\item 
If $\alpha$ is a rational number such that $0 \leq \alpha <1$ and $-(K_X+\alpha \Gamma)$ is ample, 
then it holds that  $1/3 < \alpha$.  
\end{enumerate}
\end{lem}

\begin{proof}
We fix a rational number $\alpha$ 
such that $0 \leq \alpha <1$ and $-(K_X+\alpha \Gamma)$ is ample, 
whose existence is guaranteed by Lemma \ref{l-ess-klt-case}. 

Let us show (1). 
It holds that
\[
(K_X+\Gamma) \cdot_k \Gamma \leq  (K_X+\alpha \Gamma) \cdot_k \Gamma <0, 
\]
where the first inequality follows from $\Gamma^2\leq 0$ and $0\leq \alpha <1$, 
whilst the second one holds since $-(K_X+\alpha \Gamma)$ is ample. 
Therefore, by adjunction and  Lemma \ref{l-conic}, we deduce 
$(K_X+\Gamma) \cdot_k \Gamma =\deg_k \omega_{\Gamma} = -2d_{\Gamma}$. 
Thus (1) holds. 

Let us show (2). 
For $k_{\Gamma}:=H^0(\Gamma, \MO_{\Gamma})$, 
the equation $d_{\Gamma} = [k_{\Gamma}:k]$ (Notation \ref{n-ess-klt-case}(5)) 
implies that 
\[
K_X \cdot_k \Gamma = \deg_k (\omega_X|_{\Gamma}) = d_{\Gamma} \cdot \deg_{k_{\Gamma}} (\omega_X|_{\Gamma}) \in d_{\Gamma} \Z. 
\]
Combining with  $K_X \cdot_k \Gamma>0$ (Notation \ref{n-ess-klt-case}(6)), 
we obtain $K_X \cdot_k \Gamma \geq d_{\Gamma}$. 
Hence, it holds that 
\[
(K_X+\beta \Gamma) \cdot_k \Gamma 
= (1-\beta) K_X \cdot_k \Gamma+\beta (K_X+\Gamma) \cdot_k \Gamma 
\]
\[
=(1-\beta) K_X \cdot_k \Gamma+\beta (-2d_{\Gamma})
\geq (1-\beta) d_{\Gamma} + \beta (-2d_{\Gamma})=d_{\Gamma}(1-3\beta). 
\]
Thus (2) holds.  
The assertion (3) follows from (2). 
\end{proof}

\begin{prop}\label{p-cov-deg-bound}
We use Notation \ref{n-ess-klt-case}. 
It holds that $m_{\Gamma} \leq 5.$ 
\end{prop}

\begin{proof}
We fix a rational number $\alpha$ 
such that $0 \leq \alpha <1$ and $-(K_X+\alpha \Gamma)$ is ample, 
whose existence is guaranteed by Lemma \ref{l-ess-klt-case}. 
Then the inequality $m_{\Gamma}<6$ holds by  
\[
\frac{2}{m_{\Gamma}} > \alpha > \frac{1}{3}, 
\]
where the first and second inequalities follow from  Lemma \ref{l-alpha-upper} 
and Lemma \ref{l-alpha-lower}, respectively. 
\end{proof}

To prove the main result of this subsection (Proposition \ref{p-ess-klt-bdd-torsion}), 
we first treat the case when $K(\Gamma)/K(B)$ is separable or purely inseparable.

\begin{lem}\label{l-sep-or-p-insep}
We use Notation \ref{n-ess-klt-case}. 
Let $L_B$ be a numerically trivial Cartier divisor on $B$. 
Then the following hold. 
\begin{enumerate}
\item If $K(\Gamma)/K(B)$ is a separable extension, then $\omega_B^{-1}$ is ample 
and $L_B \sim 0$. 
\item 
If $K(\Gamma)/K(B)$ is a purely inseparable morphism of degree $p^e$ 
for some $e \in \Z_{>0}$, 
then $p^e L_B \sim 0$. 
\end{enumerate}
\end{lem}

\begin{proof}
We first prove (1). 
Assume that $K(\Gamma)/K(B)$ is a separable extension. 
Let $\Gamma^N \to \Gamma$ be the normalisation of $\Gamma$. 
Set $\pi_{\Gamma^N}:\Gamma^N \to B$ to be the induced morphism. 
Since $\omega_{\Gamma}^{-1}$ is ample, so is $\omega_{\Gamma^N}^{-1}$. 
Hence we obtain $H^1(\Gamma^N, \MO_{\Gamma^N})=0$ (Lemma \ref{l-conic}). 
Thanks to the Hurwitz formula (cf. \cite[Theorem 4.16 in Section 7]{Liu02}), 
we have that $H^1(B, \MO_B)=0$, thus $\omega_B^{-1}$ is ample (Lemma \ref{l-conic}). 
In particular, the numerically trivial Cartier divisor $L_B$ is trivial, i.e. $L_B \sim 0$. 
Thus (1) holds. 

We now show (2). 
Since $K(\Gamma)/K(B)$ is a purely inseparable morphism of degree $p^e$, 
the $e$-th iterated absolute Frobenius morphism $F^e_B:B \to B$ factors through 
the induced morphism $\pi_{\Gamma^N}:\Gamma^N \to B$: 
\[
F^e_B:B \to \Gamma^N \xrightarrow{\pi_{\Gamma^N}} B. 
\]
It holds that $\pi_{\Gamma^N}^*L_B \sim 0$, hence 
$p^e L_B = (F^e_B)^* L_B \sim 0.$ Thus (2) holds. 
\end{proof}

\begin{prop}\label{p-ess-klt-bdd-torsion}
We use Notation \ref{n-ess-klt-case}. 
Let $L$ be a numerically trivial Cartier divisor on $X$. 
Then the following hold. 
\begin{enumerate}
\item 
If $p \geq 7$, then $L \sim 0$. 
\item 
If $p \in \{3, 5\}$, then $pL \sim 0$. 
\item 
If $p=2$, then $4L \sim 0.$ 
\end{enumerate}
\end{prop}

\begin{proof}
By \cite[Theorem 4.4]{Tan18a}, 
there exists a numerically trivial Cartier divisor $L_B$ on $B$ such that 
$\pi^* L_B \sim L$. 
If $K(\Gamma)/K(B)$ is separable, then Lemma \ref{l-sep-or-p-insep}(1) implies that $L \sim 0$. 
Therefore, we may assume that $K(\Gamma)/K(B)$ is not a separable extension. 
Thanks to Proposition \ref{p-cov-deg-bound}, we have 
\[
[K(\Gamma):K(B)]=m_{\Gamma} \leq 5.
\]

Let us show (1). 
Assume $p \geq 7$. 
In this case, 
there does not exist an inseparable extension $K(\Gamma)/K(B)$ with $[K(\Gamma):K(B)]\leq 5$. 
Thus (1) holds. 

Let us show (2). 
Assume $p \in \{3, 5\}$. 
Since  $K(\Gamma)/K(B)$ is not a separable extension and $[K(\Gamma):K(B)]\leq 5$, 
it holds that $K(\Gamma)/K(B)$ is a purely inseparable extension of degree $p$. 
Hence, Lemma \ref{l-sep-or-p-insep}(2) implies that $pL \sim 0$. 
Thus (2) holds. 

Let us show (3). 
Assume $p=2$. 
Since  $K(\Gamma)/K(B)$ is not a separable extension and $[K(\Gamma):K(B)]\leq 5$, 
there are the following three possibilities (i)--(iii). 
\begin{enumerate}
\item[(i)] $K(\Gamma)/K(B)$ is a purely inseparable extension of degree $2$. 
\item[(ii)] $K(\Gamma)/K(B)$ is a purely inseparable extension of degree $4$. 
\item[(iii)] $K(\Gamma)/K(B)$ is an inseparable extension of degree $4$ which is not purely inseparable. 
\end{enumerate}
If (i) or (ii) holds, then Lemma \ref{l-sep-or-p-insep}(2) implies that $4L \sim 0$. 
Hence we may assume that (iii) holds. 
Let $\Gamma^N \to \Gamma$ be the normalisation of $\Gamma$. 
Corresponding to the separable closure of $K(B)$ in $K(\Gamma)=K(\Gamma^N)$, 
we obtain the following factorisation 
\[
\Gamma^N \to B_1 \to B
\]
where $K(\Gamma^N)/K(B_1)$ is a purely inseparable extension of degree two and 
$K(B_1)/K(B)$ is a separable extension of degree two. 
In particular, $K(B_1)/K(B)$ is a Galois extension. 
Set $G:={\rm Gal} (K(B_1)/K(B))=\{{\rm id}, \sigma\}$. 
Since $L_B|_{\Gamma^N}  \sim L|_{\Gamma^N} \sim 0$ and 
the absolute Frobenius morphism $F_{B_1}:B_1 \to B_1$ factors through $\Gamma^N \to B_1$, 
it holds that $2L_B|_{B_1} \sim 0$. 
In particular, we have that $H^0(B_1, 2L_B|_{B_1}) \neq 0$. 
Fix $0 \neq s \in H^0(B_1, 2L_B|_{B_1})$. We obtain 
\[
0 \neq s \sigma(s) \in H^0(B_1, 4L_B|_{B_1})^{G}. 
\]
As $s \sigma(s)$ is $G$-invariant, 
$s \sigma(s)$ descends to $B$, i.e. there is an element 
\[
t \in H^0(B, 4L_B)
\]
such that $t|_{B_1}=s \sigma(s)$. 
In particular, we obtain $t \neq 0$, hence $4L_B \sim 0$. 
Therefore, we have $4L \sim 0$. 
\end{proof}

\subsection{General case}

We are ready to prove the main theorem of this section.

\begin{thm}\label{t-klt-bdd-torsion}
Let $k$ be a field of characteristic $p>0$. 
Let $X$ be a $k$-surface of del Pezzo type. 
Let $L$ be a numerically trivial Cartier divisor on $X$. 
Then the following hold. 
\begin{enumerate}
\item If $p \geq 7$, then $L \sim 0$. 
\item If $p \in \{3, 5\}$, then $pL \sim 0$. 
\item If $p=2$, then $4L \sim 0$. 
\end{enumerate}
\end{thm}

\begin{proof}
Replacing $k$ by $H^0(X, \MO_X)$, we may assume that $k=H^0(X, \MO_X)$. 
Furthermore, replacing $X$ by its minimal resolution, we may assume that $X$ is regular by Lemma \ref{l-dP-min-res}. 
We run a $K_X$-MMP: 
\[
\varphi:X=:X_0 \to X_1 \to \cdots \to X_n. 
\]
Since $-K_X$ is big, the end result $X_n$ is a $K_{X_n}$-Mori fibre space. 
It follows from \cite[Theorem 4.4(3)]{Tan18a} that 
there exists a Cartier divisor $L_n$ with $\varphi^*L_n \sim L$. 
Since also $X_n$ is of del Pezzo type by Lemma \ref{l-dP-under-bir-mor}, we may replace $X$ by $X_n$. 
Let $\pi:X \to B$ be the induced $K_X$-Mori fibre space. 

If $\dim B=0$, then we conclude by Theorem \ref{t-cano-bdd-torsion}.
Hence we may assume that $\dim B=1$. 
Since $X$ is a surface of del Pezzo type, 
there is an effective $\Q$-divisor such that $(X, \Delta)$ is klt and $-(K_X+\Delta)$ is ample. 
Hence any extremal ray of $\overline{\NE}(X)$ is spanned by a curve. 
Note that $\rho(X)=2$ and a fibre of $\pi:X \to B$ spans an extremal ray of $\overline{\NE}(X)$. 
Let $R=\R_{\geq 0}[\Gamma]$ be the other extremal ray, where $\Gamma$ is a curve on $X$. 
To summarise, (1)--(5) of Notation \ref{n-ess-klt-case} hold. 
There are the following three possibilities:
\begin{enumerate}
\item[(i)] $\Gamma^2 \geq 0$. 
\item[(ii)] $\Gamma^2<0$ and $K_X \cdot \Gamma \leq  0$. 
\item[(iii)] $\Gamma^2<0$ and $K_X \cdot \Gamma>0$. 
\end{enumerate}

Assume (i). 
In this case, any curve $C$ on $X$ is nef. 
Since $-(K_X+\Delta)$ is ample, also $-K_X$ is ample. 
Therefore, we conclude by Theorem \ref{t-cano-bdd-torsion}.

Assume (ii). 
In this case, $-K_X$ is nef and big. 
Again, Theorem \ref{t-cano-bdd-torsion} implies the assertion of Theorem \ref{t-klt-bdd-torsion}. 

Assume (iii). 
In this case, all the conditions (1)--(6) of Notation \ref{n-ess-klt-case} hold. 
Hence the assertion of Theorem \ref{t-klt-bdd-torsion} follows from Proposition \ref{p-ess-klt-bdd-torsion}. 
\end{proof}

\section{Results in large characteristic}

In this section, we prove the existence of geometrically normal birational models of log del Pezzo surfaces over imperfect fields of characteristic at least seven (Theorem \ref{t-dP-large-p}).
As consequences, we prove geometric integrality (Corollary \ref{c-geom-red-7}) and vanishing of irregularity for such surfaces (Theorem \ref{t-h1-vanish}).

\subsection{Analysis up to birational modification}

The purpose of this subsection is to prove Theorem \ref{t-dP-large-p}. 
To this end, we establish auxiliary results on Mori fibre spaces 
(Proposition \ref{p-dP-large-p1}, Proposition \ref{p-dP-large-p2}) 
We start by recalling the following well-known relation between the Picard rank and the anti-canonical volume of del Pezzo surfaces.

\begin{lem} \label{l-degree-picardrank}
Let $Y$ be a smooth weak del Pezzo surface over an algebraically closed field $k$. 
Then $\rho(Y) = 10 - K_Y^2.$ 
In particular, it holds that $\rho(Y) \leq 9$. 
\end{lem}
\begin{proof}
Let $Y=:Y_1 \rightarrow Y_2 \rightarrow \cdots \rightarrow Y_n=Z$ be a $K_Y$-MMP, where $Z$ is a weak del Pezzo surface endowed with a $K_Z$-Mori fibre space $Z \rightarrow B$.
It is sufficient to prove the relation $\rho(Z) = 10 - K_Z^2$, 
which is well known (cf. \cite[Theorem 1.28]{KM98}). 
\end{proof}

\begin{prop}\label{p-dP-large-p1}
Let $k$ be field of characteristic $p\geq 11$. 
Let $X$ be a regular del Pezzo $k$-surface such that 
$k=H^0(X, \MO_X)$. 
Then $X$ is smooth over $k$.  
\end{prop}

\begin{proof}
By Theorem \ref{t-p2-bound}, $X \times_k \overline{k}$ has at most canonical singularities.
By \cite[Theorem 6.1]{Sch08} such singularities are of type $A_{p^e-1}$. 
Since $X \times_k \overline{k}$ is a canonical del Pezzo surface, 
its minimal resolution $\pi \colon Y \rightarrow X \times_k \overline{k}$ is a smooth weak del Pezzo surface and we have 
\[
9 \geq \rho(Y) \geq \rho(X \times_k \overline{k})+ \sum_{x \in \text{Sing}(X \times_k \overline{k})} (p-1) \geq  \sum_{x \in \text{Sing}(X \times_k \overline{k})} 10, 
\]
where the first inequality follows from Lemma \ref{l-degree-picardrank} and the last inequality holds by $p \geq 11$. 
Thus, we obtain $\text{Sing}(X \times_k \overline{k})=\emptyset$, as desired. 
\end{proof}

\begin{prop}\label{p-dP-large-p2}
Let $k$ be field of characteristic $p>0$. 
Let $X$ be a regular $k$-surface of del Pezzo type such that $k=H^0(X, \MO_X)$. 
Assume that there is a $K_X$-Mori fibre space $\pi:X \to B$ to a 
projective regular $k$-curve $B$. 
Let $\Gamma$ be a curve which spans the extremal ray of $\overline{\NE}(X)$ 
not corresponding to $\pi$. 
Then the following hold. 
\begin{enumerate}
\item If $K_X \cdot \Gamma < 0 $ (resp. $\leq 0$), then $-K_X$ is ample (resp. nef and big). If $p \geq 5$, then $\omega_B^{-1}$ is ample and $B$ is smooth over $k$.
\item If $K_X \cdot \Gamma >0$ and $p \geq 7$, 
then $\omega_B^{-1}$ is ample and $B$ is smooth over $k$.  
\item 
If $K_X \cdot \Gamma >0$, $p \geq 7$, and $k$ is separably closed, 
then $\Gamma$ is a section of $\pi$ and $\pi$ is smooth. In particular, $X$ is smooth over $k$. 
\end{enumerate} 
\end{prop}

\begin{proof}
The first part of assertion (1) follows immediately from Kleimann's criterion for ampleness (resp. \cite[Theorem 2.2.16]{Laz04a}). 
Assume $p \geq 5$. 
The anti-canonical model $Z$ of $X$ is geometrically normal by Theorem \ref{t-p2-bound} and thus $H^1(Z, \mathcal{O}_Z)=0$.  
This implies that $H^1(X, \mathcal{O}_X)=0$ and $H^1(B, \mathcal{O}_B)=0$. 
Hence, the assertion (1) holds by Lemma \ref{l-conic} and Lemma \ref{p-Fano-curve}.

Let us show (2). 
The field extension $K(\Gamma)/K(B)$ corresponding to 
the induced morphism $\pi_{\Gamma} : \Gamma \to B$ is separable (Proposition \ref{p-cov-deg-bound}). 
Thus $B$ is a curve such that $\omega_B^{-1}$ is ample (Lemma \ref{l-sep-or-p-insep}).
Since $p>2$, $B$ is a $k$-smooth curve by Lemma \ref{p-Fano-curve}. 
Thus (2) holds. 

Let us show (3). 
It follows from Proposition \ref{p-MFS-basic}(6) that $\pi$ is a smooth morphism. 
Hence it suffices to show that $\pi_{\Gamma} : \Gamma \to B$ is a section of $\pi$. 
Since $K(\Gamma)$ is separable over $K(B)$ and $B$ is smooth over $k$, 
$K(\Gamma)$ is separable over $k$, i.e. $K(\Gamma)$ is geometrically reduced over $k$. 
Hence also $\Gamma$ is geometrically reduced over $k$. 
Since $X_{\overline k}$ is a smooth projective rational surface with $\rho(X_{\overline k})=2$, 
$X_{\overline k}$ is a Hirzebruch surface and 
$\pi_{\overline k}:X_{\overline k} \to B_{\overline k}$ is a projection. 
Since the pullback $\Gamma_{\overline k}$ of $\Gamma$ is a curve with $\Gamma_{\overline k}^2<0$ by Lemma \ref{l-ext-ray}, 
$\Gamma_{\overline k}$ is a section of $\pi_{\overline k}:X_{\overline k} \to B_{\overline k}$. 
The base change $\Gamma_{\overline k} \to B_{\overline k}$ is an isomorphism, 
hence so is the original one  $\pi_{\Gamma}:\Gamma \to B$. 
Thus (3) holds. 
\end{proof}

\begin{thm}\label{t-dP-large-p}
Let $k$ be a separably closed field of characteristic $p \geq 7$. 
Let $X$ be a $k$-surface of del Pezzo type such that $k=H^0(X, \MO_X)$. 
Then there exists a birational map $X \dashrightarrow Y$ 
to a projective normal $k$-surface $Y$ such that 
one of the following properties holds. 
\begin{enumerate}
\item $Y$ is a regular del Pezzo surface such that $k=H^0(Y, \MO_Y)$ and $\rho(Y)=1$. 
In particular, $Y$ is geometrically canonical over $k$. 
Moreover, if $p\geq 11$, then $Y$ is smooth over $k$. 
\item There is a smooth projective morphism $\pi:Y \to B$ 
such that $B \simeq \mathbb P^1_k$ and the fibre $\pi^{-1}(b)$ is isomorphic to $\mathbb P^1_{k(b)}$ 
for any closed point $b$ of $B$, 
where $k(b)$ denotes the residue field of $b$. 
In particular, $Y$ is smooth over $k$ and $Y \times_k \overline k$ is a Hirzebruch surface. 
\end{enumerate}
\end{thm}

\begin{proof}
Let $f:Z \to X$ be the minimal resolution of $X$. 
By Lemma \ref{l-dP-min-res}, $Z$ is a $k$-surface of del Pezzo type.
We run a $K_Z$-MMP: 
\[
Z=:Z_0 \to Z_1 \to \cdots \to Z_n=:Y.
\]
By Lemma \ref{l-dP-under-bir-mor}, the surfaces $Z_i$ are of del Pezzo type. 
The end result $Y$ is a $K_Y$-Mori fibre space $\pi:Y \to B$. 
If $\dim B=0$, then $Y$ is a regular del Pezzo surface, hence (1) holds by Theorem \ref{t-p2-bound} and Proposition \ref{p-dP-large-p1}. 
If $\dim B=1$, then Proposition \ref{p-dP-large-p2} implies that (2) holds. 
\end{proof}

\begin{cor}\label{c-geom-red-7}
Let $k$ be a field of characteristic $p \geq 7$. 
Let $X$ be a $k$-surface of del Pezzo type such that $k=H^0(X, \MO_X)$. 
Then $X$ is geometrically integral over $k$. 
\end{cor}

\begin{proof}
We may assume $k$ is separably closed. 
It is enough to show that $X$ is geometrically reduced \cite[Lemma 2.2]{Tan18b}. 
By Lemma \ref{l-gred-open}, we may replace $X$ by a surface birational to $X$. 
Then the assertion follows from Theorem \ref{t-dP-large-p}. 
\end{proof}

\subsection{Vanishing of $H^1(X, \MO_X)$}

In this subsection, we prove that surfaces of del Pezzo type over an imperfect field of characteristic $p \geq 7$ have vanishing irregularity.

\begin{lem}\label{l-h1-vanish}
Let $k$ be a field of characteristic $p > 0$. 
Let $X$ be a $k$-surface of del Pezzo type such that $k=H^0(X, \MO_X)$. 
If $X$ is geometrically normal over $k$, 
then it holds that $H^i(X, \MO_X)=0$ for $i>0$. 
\end{lem}

\begin{proof}
The assertion immediately follows from Lemma \ref{l-rationality}. 
\end{proof}

\begin{thm}\label{t-h1-vanish}
Let $k$ be a field of characteristic $p\geq 7$. 
Let $X$ be a $k$-surface of del Pezzo type such that $k=H^0(X, \mathcal{O}_X)$. 
Then $H^i(X, \MO_X)=0$ for $i>0$. 
\end{thm}

\begin{proof}
We may assume that $k$ is separably closed. 
Let $X \dashrightarrow Y$ be the birational morphism as in the statement 
of Theorem \ref{t-dP-large-p}. 
Lemma \ref{l-h1-vanish} implies that $H^i(Y, \MO_Y)=0$ for $i>0$. 

Let $\varphi:W \to X$ and $\psi:W \to Y$ be birational morphisms 
from a regular projective surface $W$. 
Since both $Y$ and $W$ are regular, we have that $H^i(W, \MO_W)=0$ for $i>0$. 
Then the Leray spectral sequence implies that $H^1(X, \MO_X)=0$. 
It is clear that $H^j(X, \MO_X)=0$ for $j \geq 2$. 
\end{proof}

\begin{rem}
We now give an alternative proof of Theorem \ref{t-h1-vanish}. 
We use the same notation as in \cite[Chapter 9]{FGAex}. 
Assume that $H^1(X, \MO_X) \neq 0$ and let us derive a contradiction. 
We may assume that $k$ is separably closed. 
Since $X$ is geometrically integral over $k$ (Corollary \ref{c-geom-red-7}), 
$X$ has a $k$-rational point, i.e. $X(k) \neq \emptyset$. 
By \cite[Theorem 9.2.5 and Corollary 9.4.18.3]{FGAex}, 
there exists a scheme $\mathbf{Pic}_{X/k}$ that represents 
any of the functors $\text{Pic}_{X/k}$, $\text{Pic}_{X/k, (\text{\'et})}$, and 
$\text{Pic}_{X/k, (\text{fppf})}$. 
Then $\mathbf{Pic}_{X/k}$ is a group $k$-scheme which is locally of finite type over $k$ \cite[Proposition 9.4.17]{FGAex} and its connected component $\mathbf{Pic}_{X/k}^0$ containing the identity is an open and closed group subscheme of finite type over $k$ \cite[Proposition 9.5.3]{FGAex}. 
By $H^1(X, \MO_X) \neq 0$ and $H^2(X, \MO_X)=0$, 
$\mathbf{Pic}_{X/k}^0$ is smooth and $\dim \mathbf{Pic}_{X/k}^0>0$ 
\cite[Remark 9.5.15 and Theorem 9.5.11]{FGAex}. 
Since $k$ is separably closed, $\mathbf{Pic}_{X/k}^0(k)$ is an infinite set. 
In particular, there exists a numerically trivial Cartier divisor $L$ on $X$ 
with $L \not\sim 0$. 
This contradicts Theorem \ref{t-klt-bdd-torsion}. 
\end{rem}

In characteristic zero, it is known that the image of a variety of Fano type under a surjective morphism remains of Fano type (cf. \cite[Theorem 5.12]{FG12}).
The same result is false over imperfect fields of low characteristic as shown in \cite[Theorem 1.4]{Tan}. 
We now prove that this phenomenon can appear exclusively in low characteristic.

\begin{cor}\label{c-h1-vanish}
Let $k$ be a field of characteristic $p \geq 7$. 
Let $X$ be a $k$-surface of del Pezzo type such that $k=H^0(X, \mathcal{O}_X)$ and
let $\pi \colon X \rightarrow Y$ be a projective $k$-morphism such that $\pi_* \mathcal{O}_X= \mathcal{O}_Y$. 
Then $Y$ is a $k$-variety of Fano type. 
Furthermore, if $\dim Y=1$, then $Y$ is smooth over $k$. 
\end{cor}

\begin{proof}
We distinguish two cases according to $\dim Y$.
If $\dim Y=2$, then $\pi$ is birational and we conclude by Lemma \ref{l-dP-under-bir-mor}.
If $\dim Y=1$, then thanks to the Leray spectral sequence, we have an injection: 
\[
H^1(Y, \MO_Y) \hookrightarrow H^1(X, \MO_X),
\]
where $H^1(X, \MO_X)=0$ by Theorem \ref{t-h1-vanish}. Therefore $\omega_Y^{-1}$ is ample by Lemma \ref{l-conic} and $Y$ is smooth over $k$ by Lemma \ref{p-Fano-curve}.
\end{proof}

\section{Purely inseparable points on log del Pezzo surfaces}\label{s-pi-pts}

The aim of this section is to construct purely inseparable points 
of bounded degree on log del Pezzo surfaces $X$ over $C_1$-fields of positive characteristic 
(Theorem \ref{t-ex-rat-points-dP}). 
Since we may take birational model changes, 
the problem is reduced to the case when $X$ has a Mori fibre space structure $X \to B$. 
The case when $\dim B=0$ and $\dim B=1$ are treated in 
Subsection \ref{ss1-pi-pts} and Subsection \ref{ss2-pi-pts}, respectively. 
In Subsection \ref{ss3-pi-pts}, we prove the main result 
of this section (Theorem \ref{t-ex-rat-points-dP}).

\subsection{Purely inseparable points on regular del Pezzo surfaces}\label{ss1-pi-pts}

In this subsection we prove the existence of purely inseparable points with bounded degree on geometrically normal regular del Pezzo surfaces over $C_1$-fields. 
If $K_X^2 \leq 4$,  then we apply the strategy as in \cite[Theorem IV.6.8]{Kol96} 
(Lemma \ref{l-rat-pts-low-deg}). 
We analyse the remaining cases by using a classification result given by \cite[Section 6]{Sch08} and Proposition \ref{p-insep-bdd-rat-pts}. 
We first relate the $C_r$-condition 
(for definition of $C_r$-field, see \cite[Definition IV.6.4.1]{Kol96}) 
for a field of positive characteristic to its $p$-degree.

\begin{lem}\label{l-Cr-pdeg}
Let $k$ be a field of characteristic $p>0$.
If $r$ is a positive integer and $k$ is a $C_r$-field, then $\pdeg (k) \leq r$, 
where $\pdeg (k):=\log_p [k:k^p]$. 
In particular, if $k$ is a $C_1$-field, then $\pdeg(k) \leq 1$. 
\end{lem}

\begin{proof}
Suppose by contradiction that $[k:k^p] \geq p^{r+1}$.
Let $s_1, ..., s_{p^r+1}$ be elements of $k$ which are linearly independent over $k^p$. 
Let us consider the following homogeneous polynomial of degree $p$:
\[ P:=\sum_{k=1}^{p^r+1} s_i x_i^p =s_1x_1^p + \cdots+ s_{p^r+1}x_{p^r+1}^p \in k[x_1, \dots, x_{p^r+1}]. \]
Since $s_1, ..., s_{p^r+1}$ are linearly independent over $k^p$, the polynomial $P$ has only the trivial solution in $k$.
In particular $k$ is not a $C_r$-field.
\end{proof}

We then study rational points on geometrically normal del Pezzo surfaces of degree $\leq 4$ (compare with \cite[Exercise IV.6.8.3]{Kol96}).
We need the following result. 

\begin{lem}[cf. Exercise IV.6.8.3.2 of \cite{Kol96}]\label{l-pt-dP2}
Let $k$ be a $C_1$-field. Let $S$ be a weighted hypersurface of degree 4 in $\mathbb{P}_k(1,1,1,2)$.
Then $S(k) \neq \emptyset$.
\end{lem}

\begin{proof} 
Let us recall the definition of normic forms (\cite[Definition IV.6.4.2]{Kol96}). 
A homogeneous polynomial $h \in k[y_1,  ..., y_m]$ of degree $m$ 
is called a normic form if $h=0$ has only the trivial solution in $k$. 
If $k$ has a normic form of degree two, then the same argument as in the proof of \cite[Theorem IV.6.7]{Kol96} works.

Suppose now that $k$ does not have a normic form of degree two. 
We can write $\mathbb{P}_k(1,1,1,2) = \Proj\,k[x_0, x_1, x_2, x_3]$, 
where $\deg x_0 = \deg x_1 =\deg x_2=1$ and $\deg x_3=2$. 
Let 
\[
F(x_0, x_1, x_2, x_3):= c x_3^2+f(x_0, x_1, x_2)x_3+g(x_0, x_1, x_2) \in k[x_0, x_1, x_2, x_3]
\]
be the defining polynomial of $S$, 
where $c \in k$ and $f(x_0, x_1, x_2), g(x_0, x_1, x_2) \in k[x_1, x_2, x_3]$. 
If $c=0$, then $F(0, 0, 0, 1)=0$. 
Thus, we may assume that $c \neq 0$. 
Fix $(a_0, a_1, a_2) \in k^3 \setminus \{(0, 0, 0)\}$. 
Set $\alpha:=f(a_0, a_1, a_2) \in k$ and $\beta:=g(a_0, a_1, a_2) \in k$. 
Since $h(X, Y):=cX^2 +\alpha XY +\beta Y^2$ is not a normic form, 
there is $(u, v) \in k^2 \setminus \{(0, 0)\}$ such that $h(u, v) = cu^2+\alpha uv + \beta v^2 =0$. 
Since $c \neq 0$, we obtain $v \neq 0$. 
Therefore, it holds that $F(a_0, a_1, a_2, u/v)=c(u/v)^2+\alpha (u/v)+\beta=0$, as desired. 
\end{proof}

\begin{lem}\label{l-rat-pts-low-deg}
Let $X$ be a geometrically normal regular del Pezzo surface over a $C_1$-field $k$ of characteristic $p>0$ such that $k=H^0(X, \mathcal{O}_X)$. If $K_X^2 \leq 4$, then $X(k) \neq \emptyset$.
\end{lem}

\begin{proof}
Since $X$ is geometrically normal, then it is geometrically canonical by Theorem \ref{t-classify-bc}. 
Thus we can apply Theorem \ref{t-dP-small-degree} and
we distinguish the cases according to the degree of $K_X$.

If $K_X^2=1$, then $X$ has a $k$-rational point by Proposition \ref{p-antican-ring-dP}(2).
If $K_X^2=2$, then $X$ can be embedded as a weighted hypersurface of degree 4 in $\mathbb{P}_k (1,1,1,2)$ and we apply Lemma \ref{l-pt-dP2} to conclude it has a $k$-rational point.
If $K_X ^2=3$, then $X$ is a cubic hypersurface in $\mathbb{P}^3_k$ and thus it has a $k$-rational point by definition of $C_1$-field.
If $K_X^2=4$, then $X$ is a complete intersection of two quadrics in $\mathbb{P}^4$ and thus it has a $k$-rational point by \cite[Corollary in page 376]{Lan52}.
\end{proof}

We now discuss the existence of purely inseparable points on geometrically normal regular del Pezzo surfaces over $C_1$-fields.

\begin{prop}\label{p-rat-point-p-geq-7}
Let $X$ be a regular del Pezzo surface over a $C_1$-field $k$ of characteristic $p \geq 7$ such that $k=H^0(X, \mathcal{O}_X)$. Then $X(k) \neq \emptyset$.
\end{prop}
\begin{proof}
If $X$ is a smooth del Pezzo surface, we conclude that there exists a $k$-rational point by \cite[Theorem IV.6.8]{Kol96}. 
If $p \geq 11$, then $X$ is smooth by Proposition \ref{p-dP-large-p1} and we conclude.

It suffices to treat the case when $p=7$ and $X$ is not smooth. 
By Theorem \ref{t-p2-bound}(2), $X$ is geometrically canonical.
By \cite[Theorem 6.1]{Sch08}, any singular point of 
the base change $X_{\overline k} = X \times_k \overline k$ 
is of type $A_{p^n-1}$. 
It follows from Lemma \ref{l-degree-picardrank} 
that $X_{\overline{k}}$ has a unique $A_6$ singular point.
Thus by Lemma \ref{l-degree-picardrank} we have $K_X^2  \leq 3$, 
hence Lemma \ref{l-rat-pts-low-deg} implies $X(k) \neq \emptyset$.
\end{proof}

\begin{prop} \label{p-ins-rat-point-3,5}
Let $X$ be a regular del Pezzo surface over a $C_1$-field $k$ of characteristic $p \in \left\{3,5 \right\}$ such that $k=H^0(X, \mathcal{O}_X)$. 
If $X$ is geometrically normal over $k$, then $X(k^{1/p}) \neq \emptyset$.
\end{prop}

\begin{proof}
It is sufficient to consider the case when $X$ is not smooth by \cite[Theorem IV.6.8]{Kol96}. 
By Theorem \ref{t-classify-bc}, $X_{\overline{k}}$ has canonical singularities.

If $p=5$ and $X$ is not smooth, then the singularities of $X_{\overline{k}}$ must be of type $A_4$ or $E^0_8$ according to \cite[Theorem 6.1 and Theorem 6.4]{Sch08}.
If $X_{\overline{k}}$ has one singular point of type $E_8^0$ or two singular points of type $A_4$, then $K_X^2=1$ by Lemma \ref{l-degree-picardrank}. Thus we conclude that $X$ has a $k$-rational point by Lemma \ref{l-rat-pts-low-deg}.
If $X_{\overline{k}}$ has a unique singular point of type $A_4$, it follows from Proposition \ref{p-insep-bdd-rat-pts} that $X(k^{1/p}) \neq \emptyset$.

If $p=3$ and $X$ is not smooth, then the singularities of $X_{\overline{k}}$ must be of type $A_2$, $A_8$, $E_6^0$ or $E_8^0$ according to \cite[Theorem 6.1 and Theorem 6.4]{Sch08}.
If one of the singular points is of the type $A_8$, $E_6^0$ and $E_8^0$, then $K_X^2 \leq 3$ by Lemma \ref{l-degree-picardrank} and we conclude $X(k) \neq \emptyset$ by Lemma \ref{l-rat-pts-low-deg}. 
Thus, we may assume that all the singularities of $X_{\overline k}$ are of type $A_2$. 
If there is a unique singularity of type $A_2$ on $X_{\overline k}$, 
then it follows from Proposition \ref{p-insep-bdd-rat-pts} that $X(k^{1/3}) \neq \emptyset$. 
Therefore, we may assume that 
there are at least two singularities of type $A_2$ on $X_{\overline k}$. 
Then it holds that $K_X^2 \leq 5$. 
By \cite[Table 8.5 in page 431]{Dol12}, we have that $K_X^2 \neq 5$, 
hence $K_X^2 \leq 4$. Thus Lemma \ref{l-rat-pts-low-deg} implies $X(k) \neq \emptyset$.
\end{proof}

\begin{prop} \label{p-ins-rat-point-2}
Let $X$ be a regular del Pezzo surface over a $C_1$-field $k$ of characteristic $p=2$ such that $k=H^0(X, \mathcal{O}_X)$. 
If $X$ is geometrically normal, then $X(k^{1/4}) \neq \emptyset$.
\end{prop}

\begin{proof}
It is sufficient to consider the case when $X$ is not smooth by \cite[Theorem IV.6.8]{Kol96}. 
The singularities of $X_{\overline{k}}$ are canonical by Theorem \ref{t-classify-bc}. 
Hence, by \cite[Theorem in page 57]{Sch08}, they must be of type $A_1$, $A_3$ $A_7$, $D_n^0$ with $4 \leq  n \leq 8$ or $E_n^0$ for $n=6,7,8$. 
We distinguish five cases for the singularities appearing on $X_{\overline{k}}$. 
\begin{enumerate}
\item 
There exists at least a singular point of type $A_7$, $D_n^0$ with $n \geq 5$ or $E_n^0$ for $n=6,7,8$.  
\item There are at least two singular points with one being of type $A_3$. 
\item There exists at least one singular point of type $D_4^0$. 
\item There is a unique singular point of type $A_3$. 
\item All the singular points are of type $A_1$. 
\end{enumerate}

In case (1), it holds that $K_X^2 \leq 4$. 
Hence, we obtain $X(k) \neq \emptyset$ by Lemma \ref{l-rat-pts-low-deg}. 
In case (2), 
if $K_X^2 \leq 4$, then Lemma \ref{l-rat-pts-low-deg} again implies $X(k) \neq \emptyset$. 
Hence, we may assume that $K_X^2=5$.  
Then there exist exactly two singular points $P$ and $Q$ on $X_{\overline k}$ 
such that $P$ is of type $A_3$ and $Q$ is of type $A_1$. 
However, this cannot occur by \cite[Table 8.5 at page 431]{Dol12}. 
 
In case (3) we have that $K_X^2 \leq 5$. However a $D_4^0$ singularity cannot appear on a del Pezzo of degree five according to \cite[Table 8.5 at page 431]{Dol12}. 
Thus $K_X^2 \leq 4$ and Lemma \ref{l-rat-pts-low-deg} implies $X(k) \neq \emptyset$.
In case (4), we apply Proposition \ref{p-insep-bdd-rat-pts} to conclude that $X(k^{1/4}) \neq \emptyset$.

In case (5), consider $X_{(k^{\sep})^{1/2}}$. 
By Proposition \ref{p-insep-bdd-rat-pts}, on  $X_{(k^{\sep})^{1/2}}$ there are singular points  $\left\{P_i\right\}_{i=1}^m$ of type $A_1$ such that $\kappa(P_i)=(k^{\sep})^{1/2}$ and their union $\coprod_i P_i$ is the non-smooth locus of $X_{(k^{\sep})^{1/2}}$.
Let $Y= \text{Bl}_{\coprod_i P_i} X_{(k^{\sep})^{1/2}}$ 
be the blowup of $X_{(k^{\sep})^{1/2}}$ along $\coprod_i P_i$. 
Since each $P_i$ is a $(k^{\sep})^{1/2}$-rational point 
whose base change to the algebraic closure is a canonical singularity of type $A_1$, the surface $Y$ is smooth.
Since the closed subscheme $\coprod_i P_i$ is invariant under the action of the Galois group $\text{Gal}((k^{\sep})^{1/2}/k ^{1/2})$,
the birational $(k^{\sep})^{1/2}$-morphism $Y \rightarrow X_{(k^{\sep})^{1/2}}$ descends
to a birational $k^{1/2}$-morphism $Z \rightarrow X_{k^{1/2}}$, 
where $Z$ is a smooth projective surface over $k^{1/2}$ 
whose base change to the algebraic closure is a rational surface.
It holds that $Z(k^{1/2}) \neq \emptyset$ by \cite[Theorem IV.6.8]{Kol96}, which implies $X(k^{1/2}) \neq \emptyset$. 
\end{proof}

\subsection{Purely inseparable points on Mori fibre spaces}\label{ss2-pi-pts}

In this subsection, we discuss the existence of purely inseparable points on log del Pezzo surfaces over $C_1$-fields admitting Mori fibre space structures onto curves. 
We start by recalling auxiliary results.

\begin{lem}\label{l-conic-rat-pt}
Let $k$ be a $C_1$-field and let $C$ be a regular projective curve such that $k=H^0(C, \mathcal{O}_C)$ and $-K_C$ is ample. 
Then it holds that $C \simeq \mathbb P^1_k$. 
In particular, $C(k) \neq \emptyset$.  
\end{lem}

\begin{proof}
Since $C$ is a geometrically integral conic curve in $\mathbb P^2_k$ (Lemma \ref{l-conic}), 
the assertion follows from definition of $C_1$-field. 
\end{proof}

\begin{lem} \label{l-suff-rat-pts-base}
Let $X$ be a regular projective surface over a $C_1$-field $k$ of characteristic $p>0$ such that $k=H^0(X, \mathcal{O}_X)$.
Let $\pi \colon X \rightarrow B$ be a $K_X$-Mori fibre space to a regular projective curve $B$. 
Then the following hold. 
\begin{enumerate}
\item 
Let $k \subset k'$ be an algebraic field extension. 
If $B(k') \neq \emptyset$, then $X(k') \neq \emptyset$.
\item 
If $-K_B$ is ample, then $X(k) \neq \emptyset$. 
\end{enumerate}
\end{lem}

\begin{proof}
Let us show (1). 
Let $b$ be a closed point in $B$ such that $k \subset \kappa(b) \subset k'$. 
By Proposition \ref{p-MFS-basic}, the fibre $X_b$ is a conic in $\mathbb{P}^2_{\kappa(b)}$.
By \cite[Corollary in page 377]{Lan52}, $\kappa(b)$ is a $C_1$-field, 
hence we deduce $X_{\kappa(b)}(\kappa(b)) \neq \emptyset$. 
Thus, (1) holds. 
The assertion (2) follows from Lemma \ref{l-conic-rat-pt} and (1) for the case when $k'=k$. 
\end{proof}

To discuss the case when $p=2$, 
we first handle a complicated case in characteristic two.

\begin{prop}\label{p-weird}
Let $k$ be a field of characteristic two such that $[k:k^2] \leq 2$. 
Let $X$ be a regular $k$-surface of del Pezzo type and 
let $\pi:X \to B$ be a $K_X$-Mori fibre space to a curve $B$. 
Let $\Gamma$ be a curve which spanns the $K_X$-negative extremal ray 
which is not corresponding to $\pi$. 
Assume that 
\begin{enumerate}
\item $K_X \cdot \Gamma>0$, and  
\item $K(\Gamma)/K(B)$ is an inseparable extension of degree four 
which is not purely inseparable. 
\end{enumerate}
Then $-K_B$ is ample. 
\end{prop}

\begin{proof}
We divide the proof in several steps.
\setcounter{step}{0}
\begin{step}\label{s1-weird}
In order to show the assertion of Proposition \ref{p-weird}, 
we may assume that 
\begin{enumerate}
\setcounter{enumi}{2}
\item $B$ is not smooth over $k$,
\item $\pdeg(k)=1$, i.e. $[k:k^2]=2$,  and  
\item the generic fibre of $\pi$ is not geometrically reduced. 
\end{enumerate}
\end{step}

\begin{proof}
If (3) does not hold, then $B$ is a smooth curve over $k$. 
Since $(X_{\overline k})_{\red}$ is a rational surface by Lemma \ref{l-rationality}, 
$B_{\overline k}$ is a smooth rational curve. 
Then $-K_B$ is ample, as desired. 
Thus, we may assume (3). 
From now on, we assume (3). 

If (4) does not hold, then $k$ is a perfect field. 
In this case, $B$ is smooth over $k$, which contradicts (3). 
Thus, we may assume (4). 

Let us prove the assertion of Proposition \ref{p-weird} if (5) does not hold. 
In this case, the generic fibre $X_{K(B)}$ of $\pi:X \to B$ is a geometrically integral regular conic over $K(B)$.
Thus it is smooth over $K(B)$ by Lemma \ref{p-Fano-curve}.
We use notation as in Notation \ref{n-ess-klt-case}. 
Lemma \ref{l-ess-klt-case}(8) enables us to find 
a rational number $\alpha$ such that $ 0\leq \alpha <1$ and $(X, \alpha \Gamma)$ 
is a log del Pezzo pair. 
Then Lemma \ref{l-alpha-upper}(3) implies that $\alpha m_{\Gamma} <2$. 
Since our assumption (2) implies $m_{\Gamma}=[K(\Gamma):K(B)]=4$, we have that $\alpha <1/2$. 
By the assumption (2) and $\alpha<1/2$, 
the induced pair $(X_{\overline{K(B)}}, \alpha \Gamma|_{X_{\overline{K(B)}}})$ 
on the geometric generic fibre is $F$-pure. 
It follows from \cite[Corollary 4.10]{Eji} that $-K_B$ is ample. 
Hence, we may assume that (5) holds. 
This completes the proof of Step \ref{s1-weird}. 
\end{proof}

From now on, we assume that (3)--(5) of Step \ref{s1-weird} hold. 

\begin{step}\label{s2-weird}
$X$ and $B$ are geometrically integral over $k$. 
$X$ is not geometrically normal over $k$. 
\end{step}

\begin{proof}
Since $[k:k^2]=2$, 
it follows from \cite[Theorem 2.3]{Sch10} 
that $X$ and $B$ are geometrically integral over $k$ 
(note that $\log_2 [k:k^2]$ is called the degree of imperfection for $k$ 
in \cite[Theorem 2.3]{Sch10}). 
If $X$ is geometrically normal over $k$, 
then also $B$ is geometrically normal over $k$, i.e. $B$ is smooth over $k$. 
This contradicts (3) of Step \ref{s1-weird}. 
This completes the proof of Step \ref{s2-weird}.
\end{proof}

We now introduce some notation. 
Set $k_1:=k^{1/2}$. 
By Step \ref{s2-weird}, $X \times_k k_1$ is integral and non-normal 
(cf. \cite[Proposition 2.10(3)]{Tan19}). 
Let $\nu:X_1:=(X \times_k k_1)^N \to X \times_k k_1$ be its normalisation. 
Let $X_1 \to B_1$ be the Stein factorisation of the induced morphism $X_1 \to X \to B$. 
To summarise, we have a commutative diagram 
\[
\begin{CD}
X_1 @>\nu >> X \times_k k_1 @>>> X\\
@VVV @VVV @VVV\\
B_1 @>>> B \times_k k_1 @>>> B.
\end{CD}
\]
Let $C \subset X \times_k k_1$ and $D \subset X_1$ be the closed subschemes 
defined by the conductors for $\nu$. 
For $K:=K(B)$, we apply the base change $(-) \times_B \Spec\,K$ to the above diagram: 
\[
\begin{CD}
V_1 @>>> V \times_K L  @>>> V\\
@VVV @VVV @VVV\\
\Spec\,K_1 @>>> \Spec\,L @>>> \Spec\,K, 
\end{CD}
\]
where $V:=X \times_B K$, $L:=K(B \times_k k_1) = K(B) \otimes_k k_1$, and $K_1=K(B_1)$. 
Since taking Stein factorisations commute with flat base changes, 
the morphism $V_1 \to \Spec\,K_1$ coincides the Stein factorisation of the induced morphism $V_1 \to \Spec\,K$.

\begin{step}\label{s3-weird}
$C$ dominates $B$.  
\end{step}

\begin{proof}
Assuming that $C$ does not dominate $B$, 
let us derive a contradiction. 
Since $B$ is geometrically integral over $k$ (Step \ref{s2-weird}), 
we can find a non-empty open subset $B'$ of $B$ such that $B'$ is smooth over $k$ and 
the image of $C$ on $B$ is disjoint from $B'$. 
Let $B'_1, X',$ and $X_1'$ be the inverse images of $B'$ to $B_1, X,$ and $X_1$, respectively. 
Then the resulting diagram is as follows
\[
\begin{CD}
X'_1 @>\simeq >> X' \times_k k_1 @>>> X'\\
@VVV @VVV @VV\pi' V\\
B'_1 @>\simeq >> B' \times_k k_1 @>>> B'.
\end{CD}
\]
Since $X'_1 \simeq X' \times_k k_1 = X' \times_k k^{1/2}$ is normal, 
it holds that $X'$ is geometrically normal over $k$.

Let $\pi'_{\overline k}:X'_{\overline k} \to B'_{\overline k}$ be the base change 
of $\pi'$ to the algebraic closure $\overline k$. 
Since $X'$ is geometrically normal over $k$, 
$X'_{\overline k}$ is a normal surface. 
Note that $B'_{\overline k}$ is a smooth curve. 
Since general fibres of $\pi'_{\overline k}:X'_{\overline k} \to B'_{\overline k}$
are $K_{X'_{\overline k}}$-negative and $(\pi'_{\overline k})_*\MO_{X'_{\overline k}} = \MO_{B'_{\overline k}}$, 
general fibres of $\pi'_{\overline k}$ are isomorphic to $\mathbb P^1_{\overline k}$.  
Then the generic fibre of $\pi'_{\overline k}:X'_{\overline k} \to B'_{\overline k}$ is smooth, 
hence so is the generic fibre of $\pi:X \to B$. 
This contradicts (5) of Step \ref{s1-weird}.  
This completes the proof of Step \ref{s3-weird}. 
\end{proof}

\begin{step}\label{s4-weird}
The following hold. 
\begin{enumerate}
\item[(i)] $L/K$ is a purely inseparable extension of degree two. 
\item[(ii)] $V$ is a regular conic curve on $\mathbb P^2_K$ which is not geometrically reduced over $K$. 
\item[(iii)] $V_1 \to V \times_K L$ is the normalisation of $V \times_K L$. 
\item[(iv)] $V \times_K L$ is an integral scheme which is not regular. 
\item[(v)] The restriction $D|_{V_1}$ of the conductor $D$ to $V_1$ 
satisfies $D_{V_1} =Q$, where $Q$ is a $K_1$-rational point. 
\item[(vi)] $V_1$ is isomorphic to $\mathbb P^1_{K_1}$. 
\item[(vii)] $[K_1:L]$ is a purely inseparable extension of degree two, and $K_1=K^{1/2}$. 
\end{enumerate}
\end{step}

\begin{proof}
The assertions (i)--(iii) follows from the construction. 
Step \ref{s3-weird} implies (iv). 
Let us show (v). 
For the induced morphism $\varphi:V_1 \to V$, we have that 
\[
K_{V_1}+ D|_{V_1} \sim \varphi^*K_V.
\]
Since $-K_V$ is ample, it holds that 
\[
0> \deg_{K_1}(K_{V_1}+D|_{V_1}) \geq -2+\deg_{K_1}(D|_{V_1}), 
\]
which implies $\deg_{K_1} (D|_{V_1}) \leq 1$. 
Step \ref{s3-weird} implies that $D|_{V_1} \neq 0$, 
hence $D|_{V_1}$ consists of a single rational point. 
Thus, (v) holds.

Let us show (vi). 
Since $V_1$ has a $K_1$-rational point around which $V_1$ is regular, 
$V_1$ is smooth around this point. 
In particular, Lemma \ref{l-gred-open} implies that $V_1$ is geometrically reduced. 
Then $V_1$ is a geometrically integral conic curve in $\mathbb P^2_{K_1}$. 
Therefore, $V_1$ is smooth over $K_1$.  
Since $V_1$ has a $K_1$-rational point, 
$V_1$ is isomorphic to $\mathbb P^1_{K_1}$. 
Thus, (vi) holds. 

Let us show (vii). 
The inclusion $K_1 \subset K^{1/2}$, which is equivalent to $K_1^2 \subset K$, 
follows from the fact that 
$K$ is algebraically closed in $K(V)$ and the following: 
\[
K_1^2 \subset K(V_1)^2 = K(V \times_K L)^2 = (K(V) \otimes_K L)^2 \subset K(V). 
\]

It follows from \cite[Theorem 3]{BM40} that 
the $p$-degree $\pdeg(K)$ is two, i.e. $[K^{1/2}:K]=4$ 
(note that the $p$-degree is called the degree of imperfection in \cite{BM40}). 
Hence, it is enough to show that $K_1 \neq L$. 
Assume that $K_1  = L$. 
Then $V_1$ is smooth over $L$ by (vi). 
Hence, $V \times_K L$ is geometrically integral over $L$. 
Therefore, $V$ is geometrically integral over $K$, which contradicts (5) of Step \ref{s1-weird}. 
This completes the proof of Step \ref{s4-weird}. 
\end{proof}

\begin{step}\label{s5-weird}
Set-theoretically, $C$ does not contain $\Gamma \times_k k_1$.   
\end{step}

\begin{proof}
Assuming that $C$ contains $\Gamma \times_k k_1$, 
let us derive a contradiction. 
In this case, the set-theoretic inclusion 
\[
f^{-1}(\Gamma) \subset \nu^{-1}(C) = D
\]
holds, where $f:X_1 \to X$ is the induced morphism. 
Since $B_1 \to B$ is a universal homeomorphism and 
the geometric generic fibre 
$\Gamma \times_B \Spec\,\overline{K}$ of $\Gamma \to B$ consists of two points, 
the geometric generic fibre of $D \to B_1$ contains two distinct points. 
In particular, it holds that $\deg_{K_1}(D|_{V_1}) \geq 2$. 
However, this contradicts (v) of Step \ref{s4-weird}. 
This completes the proof of Step \ref{s5-weird}. 
\end{proof}

\begin{step}\label{s6-weird}
$-K_{B_1}$ is ample. 
\end{step}

\begin{proof}
It follows from Lemma \ref{l-ess-klt-case}(8) that 
there is a rational number $\alpha$ such that $0\leq \alpha<1$ and $(X, \alpha \Gamma)$ 
is a log del Pezzo pair. 
Consider the pullback: 
\[
K_{X_1}+D +\alpha f^*\Gamma = f^*(K_X+ \alpha \Gamma). 
\]
Take the geometric generic fibre $W$ of $\pi_1:X_1 \to B_1$, 
i.e. $W=V_1 \times_{K_1} \Spec\,\overline{K}_1 \simeq \mathbb P^1_{\overline K_1}$ (Step \ref{s4-weird}(vi)). 
It is clear that $-(K_W+(D +\alpha f^*\Gamma)|_W)$ is ample. 
Since $D|_{V_1}=Q$ is a rational point (Step \ref{s4-weird}(v)), 
its pullback $D|_W=:Q_W$ to $W$ is a closed point on $W$. 
As $-(K_W+(D+ \alpha f^*\Gamma)|_W)$ is ample, 
all the coefficients of $B:=(\alpha f^*\Gamma)|_W$ must be less than one. 
Therefore, Step \ref{s5-weird} implies that $(W, (D+ \alpha f^*\Gamma)|_W)$ is $F$-pure. 
It follows from \cite[Corollary 4.10]{Eji} that 
$-K_{B_1}$ is ample. 
This completes the proof of Step \ref{s6-weird}. 
\end{proof}

\begin{step}\label{s7-weird}
$-K_B$ is ample. 
\end{step}

\begin{proof}
As $-K_{B_1}$ is ample (Step \ref{s6-weird}), 
Lemma \ref{l-conic} implies that $H^1(B_1, \MO_{B_1})=0$. 
Since $K(B_1)=K(B)^{1/2}$ (Step \ref{s4-weird}(vii)), 
the morphism $B_1 \to B$ coincides with the absolute Frobenius morphism 
of $B$. 
Hence, $B_1$ and $B$ are isomorphic as schemes. 
Thus, the vanishing $H^1(B_1, \MO_{B_1})=0$ implies $H^1(B, \MO_B)=0$. 
Then $-K_B$ is ample by Lemma \ref{l-conic}. 
This completes the proof of Step \ref{s7-weird}. 
\end{proof}
Step \ref{s7-weird} completes the proof of Proposition \ref{p-weird}. 
\end{proof}

\begin{prop}\label{p-rat-point-mfs}
Let $X$ be a regular $k$-surface of del Pezzo type over a $C_1$-field $k$ of characteristic $p>0$ such that $k=H^0(X, \mathcal{O}_X)$. 
Let $\pi \colon X \rightarrow B$ be a $K_X$-Mori fibre space 
to a regular projective curve.
Then the following hold. 
\begin{enumerate}
\item If $p \geq 7$, then $X(k) \neq \emptyset$. 
\item If $p = \left\{ 3, 5 \right\}$, then $X(k^{1/p}) \neq \emptyset$.
\item If $p = 2$, then $X(k^{1/4}) \neq \emptyset$. 
\end{enumerate}
\end{prop}

\begin{proof}
Let $R=\R_{\geq 0}[\Gamma]$ be the extremal ray of $\overline{\NE}(X)$ not corresponding to $\pi:X \to B$. 
In particular, we have $\pi(\Gamma)=B$.
We distinguish two cases:
\begin{enumerate}
\item[(I)] $K_X \cdot \Gamma \leq 0$;
\item[(II)] $K_X \cdot \Gamma >0$.
\end{enumerate}

Suppose that (I) holds. 
In this case, $-K_X$ is nef and big. 
If $p>2$, then the generic fibre $X_{K(B)}$ is a smooth conic. 
In particular, the base change $X_{\overline{K(B)}}$ is strongly $F$-regular.
By \cite[Corollary 4.10]{Eji}, $-K_B$ is ample. 
Hence, Proposition \ref{p-rat-point-mfs} implies $X(k) \neq \emptyset$. 

We now treat the case when (I) holds and $p=2$. 
Then $-K_X$ is semi-ample and big. 
Let $Z$ be its anti-canonical model. 
In particular, $Z$ is a canonical del Pezzo surface. 
By Theorem \ref{t-p2-bound}, we have $\ell_F(Z/k) \leq 2$. 
Therefore, for $k_W:=k^{1/4}$ and $W:=(Z \times_k k_W)^N_{\red}$, $W$ is geometrically normal over $k_W$. In particular, $H^0(W, \mathcal{O}_W) = k_W=k^{1/4}$.
We have the following commutative diagram 
\[
\begin{CD}
Y @> \nu >> X \times_k k^{1/4} @>>> X \\
@V f VV @VVV  @VVV\\
W @> \mu >> Z \times_{k} k^{1/4} @>>> Z \\
@VVV @VVV  @VVV\\
\Spec\,k^{1/4} @= \Spec\,k^{1/4} @>>> \Spec\,k,
\end{CD}
\]
where $\mu$ and $\nu$ are the normalisations. 
It follows from Theorem \ref{t-classify-bc} 
that $W$ is geometrically klt and $H^1(W, \mathcal{O}_W)=0$.
Since the morphism $Y \rightarrow W$ is birational and $W$ is klt by Proposition \ref{p-klt-descent}, it holds that $H^1(Y, \mathcal{O}_Y)=0$.

Consider the Stein factorisation $\pi_1 \colon Y \rightarrow B_1$ of the induced morphism 
$Y \rightarrow X \xrightarrow{\pi} B$. 
Since $H^1(Y, \mathcal{O}_Y)=0$, we conclude that $H^1(B_1, \mathcal{O}_{B_1})=0$.
In particular, since $k_W$ is a $C_1$-field, 
it holds that $B_1 \simeq \mathbb{P}^1_{k_W}$ (Lemma \ref{l-conic-rat-pt}). 
Thanks to \cite[Theorem 4.2]{Tan18b}, 
we can find an effective divisor $D$ on $Y$ such that $K_Y+D=f^*K_X.$ 
Since $-K_X$ is big, also $-K_Y$ is big. 
Fix a general $k_W$-rational point $c \in B_1$ and let $F_c$ be its $\pi_1$-fibre. 
Since we take $c$ to be general, $F_c$ avoids the non-regular points of $Y$. 
By adjunction, $\omega_{F_c}^{-1}$ is ample. 
This implies that $F$ is a conic on $\mathbb P^2_{k_W}$. 
Hence, $Y(k^{1/4})= Y(k_W) \neq \emptyset$. 
Therefore, we deduce $X(k^{1/4}) \neq \emptyset$.

We suppose (II) holds. 
We have $[K(\Gamma):K(B)] \leq 5$ by Proposition \ref{p-cov-deg-bound}. 
If $K(\Gamma)/K(B)$ is separable, then $-K_B$ is ample (Lemma \ref{l-sep-or-p-insep}). 
Then Proposition \ref{p-rat-point-mfs} implies $X(k) \neq \emptyset$. 
Hence, we may assume that $K(\Gamma)/K(B)$ is inseparable. 
If $K(\Gamma)/K(B)$ is not purely inseparable, 
then $-K_B$ is ample by Proposition \ref{p-weird}. 
Again, Proposition \ref{p-rat-point-mfs} implies $X(k) \neq \emptyset$. 
Hence, it is enough to treat the case when 
$K(\Gamma)/K(B)$ is purely inseparable. 
Since $[K(\Gamma):K(B)] \leq 5$, it suffices to prove that $X(k^{1/p^e}) \neq \emptyset$ 
for the positive integer $e$ defined by $[K(\Gamma):K(B)]=p^e$. 
Set $C:=\Gamma^N$. 
Since $\omega_{\Gamma}^{-1}$ is ample, also $-K_C$ is ample. 
Hence Proposition \ref{p-rat-point-mfs} implies $C(k') \neq \emptyset$, 
where $k':=H^0(C, \MO_C)$. 
Since 
\[
k'^{p^e} \subset K(\Gamma)^{p^e} \subset K(B), 
\]
it holds that $k'^{p^e} \subset k$. 
Therefore, we obtain $X(k^{1/p^e}) \neq \emptyset$, as desired. 
\end{proof}

\subsection{General case}\label{ss3-pi-pts}
In this subsection, 
using the results proven above, 
we prove the main result in this section (Theorem \ref{t-ex-rat-points-dP}) 
We present a generalisation of the Lang--Nishimura theorem on rational points. 
Although the argument is similar to the one in \cite[Proposition A.6]{RY00}, 
we include the proof for the sake of completeness. 
\begin{lem}[Lang-Nishimura]\label{l-inv-rat-pts-reg}
Let $k$ be a field.
Let $f : X \dashrightarrow Y$ be a rational map 
between $k$-varieties.
Suppose that $X$ is regular and $Y$ is proper over $k$. 
Fix a closed point $P$ on $X$. 
Then there exists a closed point $Q$ on $Y$ 
such that $k \subset \kappa(Q) \subset \kappa(P)$, where $\kappa(P)$ and $\kappa(Q)$ denote the residue fields.
\end{lem}

\begin{proof}
The proof is by induction on $n:=\dim X$. 
If $n=0$, then there is nothing to show. 
Suppose $n>0$. 
Consider the blowup $\pi:\text{Bl}_P X \to X$ at the closed point $P$.
Since $X$ is regular,
the $\pi$-exceptional divisor $E$ is isomorphic to $\mathbb{P}^{n-1}_{\kappa(P)}$ by \cite[Section 8, Theorem 1.19]{Liu02}. 
Consider now the induced map $f:\text{Bl}_{P} X \dashrightarrow Y$.
By the valuative criterion of properness, the map $f$ induces a rational map $E=\mathbb{P}^{n-1}_{\kappa(P)} \dashrightarrow Y$ from the $\pi$-exceptional divisor $E$.
Then by the induction hypothesis $Y$ has a closed point $Q$ whose residue field is contained in $\kappa(P)$. 
\end{proof}

\begin{thm} \label{t-ex-rat-points-dP}
Let $k$ be a $C_1$-field of characteristic $p>0$. Let $X$ be a $k$-surface of del Pezzo type such that $k=H^0(X, \mathcal{O}_X)$.
Then the following hold. 
\begin{enumerate}
\item If $p \geq 7$, then $X(k) \neq \emptyset$;
\item If $p \in \left\{ 3,5 \right\}$, then $X(k^{1/p}) \neq \emptyset$;
\item If $p =2$ , then $X(k^{1/4}) \neq \emptyset$.
\end{enumerate}
\end{thm}

\begin{proof}
Let $Y \to X$ be the minimal resolution of $X$. 
We run a $K_Y$-MMP 
$Y=:Y_0 \to Y_1 \to \cdots \to Y_n =: Z$. 
Note that the end result is a Mori fibre space. 
Thanks to Lemma \ref{l-inv-rat-pts-reg}, 
we may replace $X$ by $Z$. 
Hence it is enough to treat the following two cases. 
\begin{enumerate}
\item[(i)] $X$ is a regular del Pezzo surface with $\rho(X)=1$. 
\item[(ii)] There exists a Mori fibre space structure $ \pi \colon X \rightarrow B$ to a curve $B$.
\end{enumerate} 

Assume (i). 
By Lemma \ref{l-Cr-pdeg}, we have $\pdeg (k) \leq 1$. 
Therefore $X$ is geometrically normal by \cite[Theorem 14.1]{FS18}. 
Thus we conclude by Propositions \ref{p-rat-point-p-geq-7}, 
Proposition \ref{p-ins-rat-point-3,5}, and Proposition \ref{p-ins-rat-point-2}. 
If (ii) holds, then  
the assertion follows from Proposition \ref{p-rat-point-mfs}.
\end{proof}

\section{Pathological examples}\label{s-patho}

In this section, we collect pathological features appearing on surfaces of del Pezzo type over imperfect fields.

\subsection{Summary of known results}\label{ss1-patho}

We first summarise previously known examples of pathologies appearing on del Pezzo surfaces over imperfect fields.

\subsubsection{Geometric properties} 

We have shown that if $p \geq 7$ and $X$ is a surface of del Pezzo type, 
then $X$ is geometrically integral (Corollary \ref{c-geom-red-7}). 
We have established a partial result on geometric normality (Theorem \ref{t-dP-large-p}). 
Let us summarise known examples in small characteristic related to these properties. 
 
\begin{enumerate}
\item 
Let $\mathbb F$ be a perfect field of characteristic $p>0$ and 
let $k:=\mathbb F(t_1, t_2, t_3)$. 
Then 
\[
X:=\Proj\,k[x_0, x_1, x_2, x_3]/(x_0^p+t_1x_1^p+t_2x_2^p+t_3x_3^p)
\]
is a regular projective surface which is not geometrically reduced over $k$. 
It is easy to show that $H^0(X, \MO_X)=k$. 
If the characteristic of $k$ is two or three, then $-K_X$ is ample, 
hence $X$ is a regular del Pezzo surface. 
\item There exist a field of characteristic $p=2$ and a regular del Pezzo surface $X$ over $k$ 
such that $H^0(X, \MO_X)=k$, $X$ is geometrically reduced over $k$, and $X$ is not geometrically normal over $k$ (see \cite[Main Theorem]{Mad16}). 
\item If $k$ is an imperfect field of characteristic $p=2,3$ there exists a geometrically normal regular del Pezzo surface $X$ of Picard rank one which is not smooth (see \cite[Section 14, Equation 27]{FS18}). In \cite[Theorem 14.8]{FS18}, an example of a regular geometrically integral but geometrically non-normal del Pezzo surface of Picard rank two is constructed when $p=2$.
\item If $k$ is an imperfect field of characteristic $p \in \{2, 3\}$,  
then there exists a $k$-surface $X$ of del Pezzo type 
such that $H^0(X, \MO_X)=k$, $X$ is geometrically reduced over $k$, and $X$ is not geometrically normal over $k$ (\cite{Tan}). 
\end{enumerate}

\subsubsection{Vanishing of $H^1(X, \MO_X)$}\label{ss}

We have shown that if $X$ is a surface of del Pezzo type over a field of characteristic $p \geq 7$, 
then $H^i(X, \MO_X)=0$ for $i>0$. 
Let us summarise known examples in small characteristic which violate 
the vanishing of $H^1(X, \MO_X)$. 
\begin{enumerate}
\item 
If $k$ is an imperfect field of characteristic $p=2$, 
then there exists a regular weak del Pezzo surface $X$ such that $H^1(X, \MO_X) \neq 0$ 
(see \cite{Sch07}). 
\item 
There exist an imperfect field of characteristic $p=2$ and 
a regular del Pezzo surface $X$ such that $H^1(X, \MO_X) \neq 0$ 
(see \cite[Main theorem]{Mad16}). 
\item 
If $k$ is an imperfect field of characteristic $p \in \{2, 3\}$, 
then there exists a surface $X$ of del Pezzo type such that $H^1(X, \MO_X) \neq 0$
(see \cite{Tan}).
\end{enumerate}
\begin{rem}
Since $h^1(X, \mathcal{O}_X)$ is a birational invariant for surfaces with klt singularities,
the previous examples do not admit regular $k$-birational models which are geometrically normal.
This shows that Theorem \ref{t-dP-large-p} cannot be extended to characteristic two and three. 
\end{rem}

\subsection{Non-smooth regular log del Pezzo surfaces}

In this subsection, we construct examples of regular $k$-surfaces of del Pezzo type which are not smooth (cf. Theorem \ref{t-dP-large-p}). 

\begin{prop}\label{p-count}
Let $k$ be an imperfect field of characteristic $p>0$. Then there exists a $k$-regular surface $X$ of del Pezzo type which is not smooth over $k$.
\end{prop}

\begin{proof}
Fix a $k$-line $L$ on $\mathbb P^2_k$. 
Let $Q \in L$ be a closed point such that 
$k(Q)/k$ is a purely inseparable extension of degree $p$ 
whose existence is guaranteed by the assumption that $k$ is imperfect. 
Consider the blow-up $\pi : X \to \mathbb{P}^2_k$ at the point $Q$.
We have
\[ 
K_X = \pi^* K_{\mathbb{P}^2_k} + E\quad\text{and}\quad \widetilde{L}+E = \pi^*L,
\]
where $E$ denotes the $\pi$-exceptional divisor and $\widetilde L$ is the proper transform of $L$. Since $\widetilde L \cup E$ is simple normal crossing and the $\Q$-divisor 
\[
-(K_X+\widetilde L+\epsilon E) = \pi^*(K_X+L) - \epsilon E
\]
is ample for any $0 <\epsilon \ll 1$, 
the pair $(X, (1-\delta) \widetilde L+\epsilon E)$ is log del Pezzo for $0 \delta \ll 1$. 
Hence, $X$ is of del Pezzo type. 

It is enough to show that $X$ is not smooth. 
There exists an affine open subset $\Spec\,k[x, y] = \mathbb A^2_k$ of $\mathbb P^2_k$ 
such that $Q \in \Spec\,k[x, y]$ and the maximal ideal corresponding to $Q$ 
can be written as $(x^p-\alpha,y)$ for some $\alpha \in k \setminus k^p$. 
Let $X'$ be the inverse image of $\Spec\,k[x, y]$ by $\pi$. 
Since blowups commute with flat base changes, 
the base change $X'_{\overline k}$ is isomorphic to 
the blowup of $\Spec\,\overline{k}[x, y]$ 
along the non-reduced ideal $((x-\beta)^p, y)$, 
where $\beta \in \overline k$ with $\beta^p=\alpha$. 

After choosing appropriate coordinate, 
$X'_{\overline k}$ is isomorphic to the blowup of 
$\mathbb{A}^2_{\overline{k}}=\Spec\,\overline{k}[x', y']$ along $(x'^p, y')$. 
We can directly check that $X'_{\overline k}$ contains 
an affine open subset of the form $\Spec\,k[s, y, u]/(st-u^p)$, which is not smooth. 
\end{proof}

\begin{rem}
The surface $X$ constructed in Proposition \ref{p-count} is del Pezzo (resp. weak del Pezzo) if and only if $p=2$ (resp. $p \leq 3$).
Indeed, $-E^2=[k(Q):k]=p$ implies 
$K_X \cdot_k E = (K_X+E) \cdot_k E -E^2 = -2p+ p=-p$.
Thus the desired conclusion follows from 
\[
K_X \cdot_k \widetilde{L} = K_X \cdot_k \pi^*L - K_X \cdot_k E=-3+p.
\]
\end{rem} 
\section{Applications to del Pezzo fibrations}

In this section, we give applications of 
Theorem \ref{t-klt-bdd-torsion} and Theorem \ref{t-ex-rat-points-dP} on log del Pezzo surfaces over imperfect fields to the birational geometry of threefold fibrations.
The first application is to rational chain connectedness.

\begin{thm} \label{t-rc-3fold}
Let $k$ be an algebraically closed field of characteristic $p>0$. 
Let $\pi:V \to B$ be a projective $k$-morphism such that $\pi_*\MO_V=\MO_B$, 
$V$ is a normal threefold over $k$, and $B$ is a smooth curve over $k$. 
Assume that there exists an effective $\Q$-divisor $\Delta$ such that 
$(V, \Delta)$ is klt and $-(K_V+\Delta)$ is $\pi$-nef and $\pi$-big. 
Then the following hold. 
\begin{enumerate}
\item There exists a curve $C$ on $V$ such that $C \to B$ is surjective and 
the following properties hold. 
\begin{enumerate}
\item If $p\geq 7$, then $C \to B$ is an isomorphism. 
\item If $p \in \{3, 5\}$, then $K(C)/K(B)$ is a purely inseparable extension of degree 
$\leq p$. 
\item If $p=2$, then $K(C)/K(B)$ is a purely inseparable extension of degree 
$\leq 4$. 
\end{enumerate}
\item If $B$ is a rational curve, then $V$ is rationally chain connected. 
\end{enumerate}
\end{thm}

\begin{proof}
Let us show (1). 
Thanks to \cite[Ch. IV, Theorem 6.5]{Kol96}, $K(B)$ is a $C_1$-field. Then Theorem \ref{t-ex-rat-points-dP} implies the assertion (1). 
The assertion (2) follows from (1) and the fact that general fibres are 
rationally connected (see Lemma \ref{l-rationality}). 
\end{proof}

The second application is to Cartier divisors on Mori fibre spaces which are numerically trivial over the bases.

\begin{thm} \label{t-triv-lb-mfs-curve}
Let $k$ be an algebraically closed field of characteristic $p>0$.
Let $\pi \colon V \rightarrow B$ be a projective $k$-morphism such that $\pi_* \mathcal{O}_V = \mathcal{O}_B$, where $X$ is a $\mathbb{Q}$-factorial normal quasi-projective threefold and $B$ is a smooth curve. 
Assume there exists an effective $\mathbb{Q}$-divisor $\Delta$ such that $(V, \Delta)$ is klt and $\pi \colon V \rightarrow B$ is a $(K_V+\Delta)$-Mori fibre space.
Let $L$ be a $\pi$-numerically trivial Cartier divisor on $V$.
Then the following hold.
\begin{enumerate}
\item If $p \geq 7$, then $L \sim_{\pi} 0$.
\item If $p \in \left\{ 3, 5 \right\}$, then $p^2L \sim_{\pi} 0$.
\item If $p =2$, then $16 L \sim_{\pi} 0$.
\end{enumerate}
\end{thm}

\begin{proof}
We only prove the theorem in the case when $p=2$, since the other cases are similar and easier.
Since the generic fibre $V_{K(B)}$ is a $K(B)$-surface of del Pezzo type, we have by Theorem \ref{t-klt-bdd-torsion} that $4L|_{V_{K(B)}} \sim 0$.
Therefore, $4L$ is linearly equivalent to a vertical divisor, 
i.e. we have 
\[
4L\sim \sum_{i=1}^r \ell_i D_i,
\]
where $\ell_i \in \Z$ and $D_i$ is a prime divisor such that $\pi(D_i)$ is a closed point $b_i$.

Since $\rho(V/B)=1$ and $V$ is $\Q$-factorial, 
all the fibres of $\pi$ are irreducible. 
Hence, we can write $\pi^*(b_i)=n_i D_i$ for some $n_i \in \Z_{>0}$. 
Let $m_i$ be the Cartier index of $D_i$, 
i.e. the minimum positive integer $m$ such that $mD_i$ is Cartier. 
Since the divisor $\pi^*(b_i)=n_i D_i$ is Cartier, then there exists $r_i \in \Z_{>0}$ such that $n_i = r_i m_i$. 

We now prove that $r_i$ is a divisor of $4$.
Since $K(B)$ is a $C_1$-field and the generic fibre is a surface of del Pezzo type, 
we conclude by Theorem \ref{t-ex-rat-points-dP} that there exists a curve $\Gamma$ on $V$ such that the degree $d$ of the morphism $\Gamma \rightarrow B$ is a divisor of 4.
By the equation 
\[
r_i \cdot (m_iD_i) \cdot \Gamma = n_iD_i \cdot \Gamma = \pi^*(b_i) \cdot \Gamma = d, 
\]
$r_i$ is a divisor of $4$. 

Therefore, it holds that $4m_iD_i \sim_{\pi} 0$. 
On the other hand, the divisor $4L=\sum_{i=1}^r \ell_i D_i$ is Cartier, 
hence we have that $\ell_i = s_im_i$ for some $s_i \in \Z$. 
Therefore it holds that 
\[
16L\sim \sum_{i=1}^r 4\ell_i D_i\sim  \sum_{i=1}^r s_i (4 m_i D_i) \sim_{\pi} 0,
\]
as desired. 
\end{proof}

\end{document}